\documentclass[12pt, reqno]{amsart}
\usepackage{amsmath, amstext, amsbsy, amssymb, amscd}
\usepackage{amsthm, latexsym, graphics, float}

\newtheorem{theorem}{Theorem}[section]
\newtheorem{lemma}[theorem]{Lemma}
\newtheorem{proposition}[theorem]{Proposition}
\newtheorem{corollary}[theorem]{Corollary}
\theoremstyle{definition}
\newtheorem{definition}[theorem]{Definition}
\newtheorem{example}[theorem]{Example}

\theoremstyle{remark}
\newtheorem{remark}[theorem]{Remark}

\numberwithin{equation}{section}

\newcommand{\C}{ \mathbb C }

\newcommand{\End}{{\rm End}}
\newcommand{\fock}{{\mathbb H}_{\mathbb T}}

\newcommand{\hil}[1]{{S^{[#1]}}}

\newcommand{\Hom}{{\rm Hom}}

\newcommand{\la}{\lambda}

\newcommand{\Supp}{{\rm Supp}}

\newcommand{\vac}{|0\rangle}
\newcommand{\w}{\tilde}
\newcommand{\wa}{\tilde {\mathfrak a}}
\newcommand{\wt}{\tilde {\mathfrak t}}
\newcommand{\W}{\widetilde}
\newcommand{\Wq}{\widetilde Q}
\newcommand{\Wh}{\widetilde {\mathbb H}}
\newcommand{\Wfock}{\widetilde {\mathbb H}_S}
\newcommand{\Wft}{\widetilde {\mathbb H}_{\mathbb T}}
\newcommand{\Sn}{S^{[n]}}

\newcommand{\Z}{ \mathbb Z }

{\vskip-\lastskip\medskip
  \noindent
  {\em #1.}\enspace
  }%
{\qed\par\medskip
  }

\begin{document}
\title[Incidence Hilbert schemes and loop algebras]
      {Equivariant cohomology of incidence Hilbert schemes
       and loop algebras}

\author[Wei-Ping Li]{Wei-Ping Li$^1$}
\address{Department of Mathematics, HKUST, Clear Water Bay, Kowloon,
Hong Kong} \email{mawpli@ust.hk}
\thanks{${}^1$Partially supported by the grant CERG601905}

\author[Zhenbo Qin]{Zhenbo Qin$^2$}
\address{Department of Mathematics, University of Missouri, Columbia,
MO 65211, USA} \email{zq@math.missouri.edu}
\thanks{${}^2$Partially supported by an NSF grant}


\begin{abstract}
Let $S$ be the affine plane $\C^2$ together with an appropriate
$\mathbb T = \C^*$ action. Let $\hil{m,m+1}$ be the incidence
Hilbert scheme. Parallel to \cite{LQ}, we construct an infinite
dimensional Lie algebra that acts on the direct sum
$$\Wft = \bigoplus_{m=0}^{+\infty}H^{2(m+1)}_{\mathbb T}(S^{[m,m+1]})$$
of the middle-degree equivariant cohomology group of $\hil{m,m+1}$. 
The algebra is related to the loop algebra of an infinite dimensional
Heisenberg algebra. In addition, we study the transformations 
among three different linear bases of $\Wft$. 
Our results are applied to the ring structure of the
ordinary cohomology of $\hil{m,m+1}$ and to the ring of symmetric
functions in infinitely many variables.
\end{abstract}

\maketitle
\date{}

\section{\bf Introduction}

Let $S$ be the affine plane $\C^2$ together with the $\mathbb T = \C^*$ 
action
\begin{eqnarray*}
a(w, z) = (a w, a^{-1}z), \qquad\quad  a \in {\mathbb T}
\end{eqnarray*}
on the coordinate functions $w$ and $z$ of $S$. This $\mathbb T$-action
on $S$ induces a $\mathbb T$-action on the Hilbert scheme $S^{[n]}$
of $n$-points on $S$. The $\mathbb T$-fixed points in $S^{[n]}$
are of the form $\xi_\la$ where $\la$ denotes partitions of $n$.
In \cite{Na2, Na3, Vas, LQW1, LQW2}, the equivariant cohomology 
$H^*_{\mathbb T}(S^{[n]})$ of the Hilbert scheme $S^{[n]}$
has been studied via representation theory. 
A generalization of Nakajima's work \cite{Na1} to the equivariant 
cohomology $H^*_{\mathbb T}(\hil{n})$ shows in \cite{Vas} that 
the space
\begin{eqnarray*}
\mathbb H_{\mathbb T}=\bigoplus_{n=0}^{+\infty}H^{2n}_{\mathbb
T}(\hil{n})
\end{eqnarray*}
is an irreducible representation of a Heisenberg algebra generated
by the linear operators $\mathfrak a_n^{\mathbb T}$, $n\in \mathbb
Z$ in $\text{End}(\mathbb H_{\mathbb T})$. As a consequence, it
induces a linear isomorphism
\begin{eqnarray*}
\Phi: \,\, \mathbb H_{\mathbb T} \to \Lambda \otimes_\Z \C
\end{eqnarray*}
where $\Lambda$ is the vector space of symmetric functions in
infinitely many variables (see p.19 of \cite{Mac}). More
specifically, let $C=C^z$ be the $z$-axis of $S$. 
The homomorphism $\Phi$ maps $\mathfrak a_{-\la}^{\mathbb T}\vac$ 
and $[L^\la C]$ (defined in (\ref{a_la}) and (\ref{LlaC}))
to the power-sum symmetric function $p_\la$ and
the monomial symmetric function $m_\la$ respectively.
Here $[\cdot]$ denotes the equivariant fundamental cohomology class. 

A new feature of the equivariant setup is the existence of 
the $\mathbb T$-fixed points. By the localization theorem,
the ring structure of $H^*_{\mathbb T}(\hil{n})$ is easy to
describe using the fixed points $\xi_\la \in \hil{n}$. 
Note that $\Lambda$ is a ring as well by the {\it usual} 
multiplication of functions. However, $\Phi$ is not a ring 
isomorphism. On the other hand, if we define a new ring structure on
$\Lambda$ by requiring $s_\la\cdot s_\mu=\delta_{\la, \mu}h(\la)s_\la$
for the Schur functions $s_\la$ and $s_\mu$, then $\Phi$ is a ring
isomorphism from $\mathbb H_{\mathbb T}$ to $(\Lambda, \cdot)$. 
Here $h(\la)$ denotes the hook number of the Young diagram 
associated to $\la$. In fact, 
$\Phi$ maps the fixed point class $(-1)^{|\la|}/h(\la)\cdot [\la]$
to the Schur function $s_\la$. This is an extra property gained by
going to equivariant cohomology (see \cite{Vas}).

In this paper, we study the equivariant cohomology 
$H^*_{\mathbb T}(S^{[n, n+1]})$ of the incidence Hilbert scheme 
$S^{[n, n+1]}$ which is defined by
\begin{eqnarray*}
\hil{n, n+1}=\{ (\xi,\xi')  \,|\, \xi \subset \xi' \} \subset
\hil{n}\times \hil{n+1}.
\end{eqnarray*}
It is known from \cite{Ch1, Tik} that the incidence Hilbert scheme
$\hil{n, n+1}$ is irreducible, smooth and of dimension $2(n+1)$. 
Following \cite{LQ}, 
we construct the Heisenberg operators $\wa_n^{\mathbb T}, 
n \in \mathbb Z$ and the translation operator $\wt^{\mathbb T}$
on the space
\begin{eqnarray*}      
\Wft = \bigoplus_{n=0}^{+\infty} H_{\mathbb T}^{2(n+1)}(S^{[n, n+1]}).
\end{eqnarray*}
Let ${\w {\mathfrak h}}_{\mathbb T}$ be the Heisenberg algebra 
generated by the operators $\wa_n^{\mathbb T}$, $n \in \Z$.
The loop algebra of ${\w {\mathfrak h}}_{\mathbb T}$ is the space
$\C [ u, u^{-1}] \otimes_\C {\w {\mathfrak h}}_{\mathbb T}$ 
together with the Lie bracket
\begin{eqnarray*}  
[u^m \otimes g_1,  u^n \otimes g_2] =  u^{m+n} \otimes [g_1, g_2].
\end{eqnarray*}

\begin{theorem}  \label{intro_thm1}
The space $\Wft$ is a representation of the Lie
algebra $\C[u^{-1}]\otimes_\C \w{\mathfrak h}_{\mathbb T}$ with 
a highest weight vector being the vacuum vector
\begin{eqnarray*}  
\vac = [C^z] \in H_{\mathbb T}^2(\hil{0,1}) 
= H_{\mathbb T}^2(S) = H_{\mathbb T}^2(\C^2)
\end{eqnarray*}
where $u^{-1}$ acts via $\wt^{\mathbb T}$, and $C^z$ denotes 
the $z$-axis of $S = \C^2$.
\end{theorem}

It follows that a linear basis of the space $\Wft$ is given by
\begin{eqnarray*}   
  \W {\mathcal B}_2
= \left \{ \big (\wt^{\mathbb T} \big )^{i} \,
   \wa_{-\nu}^{\mathbb T} \vac \right \}_{i \ge 0, \, \nu}.
\end{eqnarray*}
On the other hand, the $\mathbb T$-fixed points of 
$\hil{n, n+1}$ are of the form $\xi_{\la, \mu} = 
(\xi_\la, \xi_\mu)$ where $\la$ and $\mu$ denote partitions 
of $n$ and $(n+1)$ respectively, and the Young diagram of 
$\la$ is contained in the Young diagram of $\mu$. Such a 
pair $(\la, \mu)$ of partitions is defined to be 
an {\it incidence} pair. For an incidence pair $(\la, \mu)$,
let
\begin{eqnarray*}      
[\la, \mu] = t^{-(n+1)} \cup [\xi_{\la, \mu}]
\in H_{\mathbb T}^{2(n+1)}(S^{[n, n+1]})
\end{eqnarray*}
where $t$ is the character associated to the $1$-dimensional
standard module $\theta$ of $\mathbb T$ on which $a \in \mathbb T$ 
acts as multiplication by $a$. By the localization
theorem, the ring structure of $H_{\mathbb T}^{2(n+1)}(S^{[n, n+1]})$
(and hence of $\Wft$) can be easily described in terms of the
classes $[\la, \mu]$. In addition, these classes form
another linear basis of $\Wh_{\mathbb T}$:
\begin{eqnarray*}   
  \W {\mathcal B}_1
= \big \{ [\la, \mu] \big \}_{(\la, \mu) \,\, \text{\rm incidence}}.
\end{eqnarray*}

\begin{theorem}  \label{intro_thm2}
There exists an algorithm to express each element $\big
(\wt^{\mathbb T} \big )^{i} \, \wa_{-\nu}^{\mathbb T} \vac$ in the
linear basis $\W {\mathcal B}_2$ as a linear combination of the
elements in the linear basis $\W {\mathcal B}_1$.
\end{theorem}

This theorem implies that the ring structure of $\Wft$ can 
also be described (implicitly) in terms of the elements
in the linear basis $\W {\mathcal B}_2$. The main idea in
proving Theorem~\ref{intro_thm2} is to introduce a third 
linear basis of the space $\Wh_{\mathbb T}$:
\begin{eqnarray*}   
\W {\mathcal B}_3 
= \left \{ [\W L^{\la, \mu} C] \right \}
\end{eqnarray*}
where $\W L^{\la, \mu} C$ is defined by (\ref{LlamuC}).
We show that there exist algorithms to express every element 
in $\W {\mathcal B}_2$ as a linear combination of the
elements in $\W {\mathcal B}_3$ and to express every element 
in $\W {\mathcal B}_3$ as a linear combination of the
elements in $\W {\mathcal B}_1$.

There are two applications of our results. The first is
to describe the ordinary cohomology ring $H^*(S^{[n, n+1]})$ 
of the incidence Hilbert scheme $S^{[n, n+1]}$. The second
is to the ring of symmetric functions. Indeed,
define a linear isomorphism
\begin{eqnarray*}  
\W \Phi: \Wh_{\mathbb T} \to \Lambda \otimes_\Z \C[v]
\end{eqnarray*}
by sending $\big (\wt^{\mathbb T} \big )^{i} \, \wa_{-\la}^{\mathbb
T} \vac$ to $p_\la \otimes v^i$. Then the ring structure on 
$\Wh_{\mathbb T}$ induces a ring structure on 
$\Lambda \otimes_\Z \C[v]$ such that $\Lambda \otimes_\Z \C 
\subset \Lambda \otimes_\Z \C[v]$ becomes a subring of
$\Lambda \otimes_\Z \C[v]$. Moreover, we have a commutative
diagram of ring homomorphisms:
\begin{eqnarray*}
\CD \mathbb H_{\mathbb T} @>{\Phi}>>
    \Lambda\otimes_{\mathbb Z}\mathbb C  \\
@VVV @VV{\iota}V \\
\Wh_{\mathbb T} @>{\W \Phi}>>\Lambda \otimes_\Z \C[v]
\endCD
\end{eqnarray*}
respecting the Heisenberg algebra actions on $\mathbb H_{\mathbb T}$ 
and $\Wh_{\mathbb T}$, where $\iota$ denotes the inclusion map.
It is natural for us to ask what the induced ring structure on
$\Lambda \otimes_\Z \C[v]$ is in the realm of symmetric functions.

The paper is organized as follows. In \S 2, we study the 
equivariant aspects of the incidence Hilbert scheme $S^{[n, n+1]}$,
including a description of the $\mathbb T$-fixed points, 
the generating function for the Betti numbers, 
a $\mathbb T$-invariant cell decomposition, 
the equivariant Zariski tangent spaces at the fixed points, 
and a bilinear pairing. In \S 3, we construct the loop algebra
action on the space $\Wft$, and compare it with the Heisenberg
algebra action on the space $\fock$. In \S 4, we study the
transformations among the three linear bases $\W {\mathcal B}_1$,
$\W {\mathcal B}_2$ and $\W {\mathcal B}_3$ of $\Wft$.
In \S 5, the two applications mentioned above are addressed.
In \S 6 (the Appendix), we prove Lemma~\ref{lemma_eT}.

\medskip\noindent
{\bf Conventions.} We use $\la$ and $\mu$ to denote partitions of
$n$ and $(n+1)$ respectively. The sign $\hbox{ }\widetilde{}\hbox{
}$, in the case of  cohomology and operators, is  for  the incidence
Hilbert schemes $S^{[n,n+1]}$. The sign $\hbox{ }{}^\prime\hbox{ }$,
in the case of equivaiant cohomology,  is  for the localized
equivariant cohomology.

\bigskip\noindent
{\bf Acknowledgments.} The first author thank the
Department of Mathematics at the University of Missouri for the
Miller Scholarship which made his visit there in February and June
of 2006 possible and MSRI at Berkeley for its support.

\section{\bf The equivariant setup for incidence Hilbert schemes}
\label{sec:setup}

When a smooth algebraic variety $X$ admits a torus $\mathbb C^*$
action, one can study its equivariant cohomology $H^*_{\mathbb
C^*}(X)$. It is known that the localized $H^*_{\mathbb
C^*}(X)^\prime$ has extra properties coming from the fixed points.
In the case of Hilbert scheme $\hil{n}$ of points  on a surface $S$,
this provides a much richer structure on the equivariant cohomology
of $\hil{n}$ than the ordinary cohomology \cite{Vas, Na3,LQW1, LQW2}.

Besides the Hilbert scheme of points $\hil{n}$ for a surface $S$, 
the incidence Hilbert scheme $\hil{n,n+1}$ for a surface $S$
is the only class of (generalized or nested) Hilbert schemes of 
points on smooth varieties of dimension bigger than one 
which are smooth for all $n$ (see \cite{Ch1}). It has a
nice generating function of Betti numbers. When the surface $S$
is $\mathbb C^2$, the torus $\mathbb C^*$ action on $\hil{n,n+1}$
was studied in details in \cite{Ch1}. In this section, we follow
Cheah's approach to the equivariant tangent spaces of the fixed points.
We calculate the generating function of the Betti numbers of
$\hil{n, n+1}$, study in details the equivariant tangent spaces, and
hence determine the ring structure of the localized equivariant
cohomology $H^*_{\mathbb C^*}(\hil{n, n+1})^\prime$ in terms of 
the fixed points. It turns out that it is more natural to 
work on a modified cohomology ring, as illustrated in \cite{Vas}, 
which will be the material in the last subsection. 
We draw a special attention to three different torus actions 
on $S = \mathbb C^2$ in (\ref{T_act}), (\ref{T+_act}) and 
(\ref{T-_act}) which serve for different purposes.

\subsection{The equivariant homology and cohomology}
\label{subsect_equi} $\,$
\par

Let $\mathbb T = \C^*$, and let $\theta$ be the $1$-dimensional
standard module of $\mathbb T$ on which $a \in \mathbb T$ acts as
multiplication by $a$, and let $t$ be the associated character. Then
the representation ring $\mathcal R(\mathbb T)$ is isomorphic to
$\mathbb Z[t, t^{-1}]$.

Let $X$ be an algebraic variety acted  by $\mathbb T$. Let
$H^*_{\mathbb T}(X)$ and $H^{\mathbb T}_*(X)$ be the equivariant
cohomology and the equivariant homology with $\C$-coefficient
respectively. Note that $H^*_{\mathbb T}(pt)=H^*(B{\mathbb
T})=\mathbb C[t]$. Then there exist bilinear maps
\begin{eqnarray*}
&\cup:& H^i_{\mathbb T}(X) \otimes H^j_{\mathbb T}(X)
\to H^{i+j}_{\mathbb T}(X), \\
&\cap:& H^i_{\mathbb T}(X) \otimes H^{\mathbb T}_j(X) \to
H_{i+j}^{\mathbb T}(X).
\end{eqnarray*}
If $X$ is of pure dimension, then there exists a linear map
\begin{eqnarray*}
D: H^i_{\mathbb T}(X) \to H^{\mathbb T}_i(X).
\end{eqnarray*}
If $X$ is smooth of pure dimension, then $D$ is an isomorphism. When
$f: Y \to X$ is a ${\mathbb T}$-equivariant and proper morphism of
varieties, we have a Gysin homomorphism
\begin{eqnarray*}
f_!: H_*^{\mathbb T}(Y) \to H_*^{\mathbb T}(X)
\end{eqnarray*}
of equivariant homology. Moreover, when both $Y$ and $X$ are smooth
of pure dimension, we have the Gysin homomorphism
\begin{eqnarray*}
D^{-1} \, f_! \, D: H^*_{\mathbb T}(Y) \to H^*_{\mathbb T}(X)
\end{eqnarray*}
of equivariant cohomology, which will still be denoted by $f_!$.

\subsection{\bf Incidence Hilbert schemes of points on surfaces}
\label{subsect_inc} $\,$
\par

Let $S$ be a smooth complex surface, and $\Sn$ be the Hilbert scheme
of points in $S$. An element in $\Sn$ is represented by a length-$n$
$0$-dimensional closed subscheme $\xi$ of $S$. For $\xi \in \Sn$,
let $I_{\xi}$ be the corresponding sheaf of ideals. It is well known
that $\Sn$ is a nonsingular complex variety of dimension $2n$.
Sending an element in $\Sn$ to its support in the symmetric product
${\rm Sym}^n(S)$, we obtain the Hilbert-Chow morphism $\pi_n: \Sn
\rightarrow {\rm Sym}^n(S)$, which is a resolution of singularities.
Let
\begin{eqnarray*}
\mathcal Z_n =\{(\xi, s) \in S^{[n]} \times S \mid s \in
\Supp(\xi) \}
\end{eqnarray*}
be the universal codimension-$2$ subscheme in $S^{[n]} \times S$.

Fix a point $s \in S$. For $m \ge 0$ and $n > 0$, we define two
closed subsets:
\begin{eqnarray}
   M_m(s)
&=&\{ \xi \in \hil{m} |\, \Supp(\xi) = \{s\}\},
               \label{mms}  \\
   M_{m, m+n}(s)
&=&\{ (\xi, \xi') |\, \xi \subset \xi' \}
   \subset M_m(s) \times M_{m+n}(s).  \label{mmns}
\end{eqnarray}
It is known that $M_{m, m+1}(s)$ and $M_{m+1}(s)$ are irreducible
with
\begin{eqnarray}  \label{dim_mms}
\dim M_{m, m+1}(s) = \dim M_{m+1}(s) = m.
\end{eqnarray}

The incidence Hilbert scheme $\hil{n, n+1}$ is defined by
\begin{eqnarray}  \label{Hnn1}
\hil{n, n+1}=\{ (\xi,\xi')  \,|\, \xi \subset \xi' \} \subset
\hil{n}\times \hil{n+1}.
\end{eqnarray}
It is known from \cite{Ch1, Tik} that the incidence Hilbert scheme
$\hil{n, n+1}$ is irreducible, smooth and of dimension $2(n+1)$. 
In fact, we have
\begin{eqnarray}   \label{nn+1cong}
\hil{n, n+1} \cong \widetilde{S^{[n]} \times S}
\end{eqnarray}
where $\widetilde{S^{[n]} \times S}$ denotes the blowup of $\hil{n}
\times S$ along the subscheme ${\mathcal Z}_n$ (see \cite{ES2}).
Note that sending a pair $(\xi,\xi') \in \hil{n, n+1}$ to the
support of $I_\xi/I_{\xi'}$ yields a morphism:
\begin{eqnarray}   \label{rho_n}
\rho_n: \hil{n, n+1} \to S
\end{eqnarray}
which is also the composition of the isomorphism (\ref{nn+1cong})
and the projection
\begin{eqnarray*}
\widetilde{S^{[n]} \times S} \to S^{[n]} \times S \to S.
\end{eqnarray*}

\subsection{The torus action on the incidence Hilbert schemes}
\label{subsect_torus} $\,$
\par

Let $S=\mathbb C^2$. Then the $2$-dimensional complex torus
${\mathbb T}^2 = (\mathbb C^*)^2$ acts on the affine coordinate
functions $w$ and $z$ of $S$ by
\begin{eqnarray}   \label{T2_act}
(a, b)w=a w, \quad (a, b)z=bz\qquad (a, b) \in {\mathbb T}^2.
\end{eqnarray}
It induces ${\mathbb T}^2$-actions on both $S^{[n]}$ and $S^{[n,
n+1]}$. It is known from \cite{ES1} that the ${\mathbb T}^2$-fixed
points in $S^{[n]}$ are parametrized by the partitions of $n$. Let
$\la$ be a partition of $n$ (denoted by $\la \vdash n$), and let
$\xi_{\la}$ be the ${\mathbb T}^2$-fixed point on $S^{[n]}$
corresponding to $\la$. If $\la = (\la_1 \ge \la_2 \ge \ldots \ge
\la_r)$ with $\la_1 + \ldots + \la_r = n$, then we have
\begin{eqnarray}   \label{Ixila}
I_{\xi_{\la}} = (w^{\la_1}, zw^{\la_2}, \ldots, z^{r-1}w^{\la_r},
z^r).
\end{eqnarray}

The multiplicity of part $i$ in a partition $\mu$ is denoted by
$m_i(\mu)$, or simply by $m_i$ if there is no confusion. Using these
multiplicities, we can also express $\mu$ as:
\begin{eqnarray*}
\mu = (1^{m_1(\mu)}2^{m_2(\mu)} \cdots i^{m_i(\mu)} \cdots ) =
(1^{m_1}2^{m_2} \cdots i^{m_i} \cdots ).
\end{eqnarray*}

\begin{lemma}  \label{fp_inc}
The ${\mathbb T}^2$-fixed points in $S^{[n, n+1]}$ are of the form
$(\xi_\la, \xi_\mu)$ where
\begin{eqnarray}      \label{fp_inc.1}
   \mu
&=&\big ( \cdots (i-1)^{m_{i-1}} i^{m_i}
       (i+1)^{m_{i+1}}  \cdots \big )  \,\, \vdash \,\, (n+1)
       \nonumber           \\
   \la
&=&\mu^{(i)} := \big ( \cdots (i-2)^{m_{i-2}}(i-1)^{m_{i-1}+1}
       i^{m_i-1} (i+1)^{m_{i+1}} \cdots \big )
\end{eqnarray}
for some $i \ge 1$ with $m_i > 0$ (the parts $(i-2)^{m_{i-2}},
(i-1)^{m_{i-1}}$ and $(i-1)^{m_{i-1}+1}$ do not appear if $i=1$).
\end{lemma}
\begin{proof}
Since the ${\mathbb T}^2$-fixed points in $S^{[n]}$ are of the form
$\xi_\la$ with $\la \vdash n$, the ${\mathbb T}^2$-fixed points in
$S^{[n, n+1]}$ are of the form $(\xi_\la, \xi_\mu)$ where $\la
\vdash n$, $\mu \vdash (n+1)$ and $I_{\xi_{\mu}} \subset
I_{\xi_{\la}}$. Let $\la = (\la_1 \ge \la_2 \ge \ldots)$ and $\mu =
(\mu_1 \ge \mu_2 \ge \ldots)$. By (\ref{Ixila}),
\begin{eqnarray*}
\mu_1 \ge \la_1, \,\, \mu_2 \ge \la_2, \ldots .
\end{eqnarray*}
Since $\sum_j \la_j = n$ and $\sum_j \mu_j = n+1$, there exists
$j_0$ satisfying
$
\mu_{j_0} = \la_{j_0} +1
$
and $\mu_j = \la_j$ whenever $j \ne j_0$. This is equivalent to
(\ref{fp_inc.1}).
\end{proof}

\begin{definition}  \label{step}
(i) The {\it step length} $s(\mu)$ of a partition $\mu$ is defined
to be
\begin{eqnarray*}
s(\mu) \,\, = \,\, \#\{i|\,\, m_i(\mu) > 0\};
\end{eqnarray*}

(ii) If $(\xi_\la, \xi_\mu) \in \hil{n, n+1}$, then $(\la, \mu)$ is
defined to be an {\it incidence pair}. Put
\begin{eqnarray*}
\xi_{\la, \mu} = (\xi_\la, \xi_\mu).
\end{eqnarray*}
\end{definition}

Fix $\mu \vdash (n+1)$. By Lemma~\ref{fp_inc}, the number of
${\mathbb T}^2$-fixed points in $S^{[n, n+1]}$ which are of the form
$(\xi_\la, \xi_\mu)$ is precisely equal to the step length of the
partition $\mu$.

Let ${\mathbb T}_+ = {\mathbb T}_- = \C^*$. Consider the three
actions on coordinate functions of $S$:
\begin{eqnarray}   
a(w, z) &=&(a w, a^{-1}z), \qquad\quad  a \in {\mathbb T}, 
                        \label{T_act}   \\
a(w, z) &=&(a^{u} w, a^{v}z),
            \qquad \quad a \in {\mathbb T}_+, 
                         \label{T+_act} \\
a(w, z) &=&(a^{-u} w, a^{-v}z), \qquad a \in {\mathbb T}_-    
                         \label{T-_act}
\end{eqnarray}
where $0 < u \ll v$. We regard them as three $1$-dimensional
subgroups of ${\mathbb T}^2$. The fixed points in $S^{[n, n+1]}$
under the action of ${\mathbb T}$ (respectively, ${\mathbb T}_+$ 
and ${\mathbb T}_-$) are exactly the same as those given by
Lemma~\ref{fp_inc}.
\subsection{\bf The generating function for the Betti numbers}
\label{subsect_betti} $\,$
\par

When the surface $S$ is projective, the generating function for the
Betti numbers of the incidence Hilbert schemes $\hil{n, n+1}$ has
been determined by Cheah \cite{Ch2}:
\begin{eqnarray}  \label{Hnn1_bi}
& &\sum_{n=0}^{+\infty} \left (\sum_{i} (-1)^i
      b_i\big ( \hil{n, n+1} \big ) z^i\right )q^n    \nonumber  \\
&=&\left (\sum_{i} (-1)^i b_i(S) z^i \right ) \cdot
   \frac{1}{1-z^2q} \cdot  \prod_{n=1}^{+\infty}
   \prod_{i} \left (\frac{1}{1-z^{2n-2+i}q^n}
   \right )^{(-1)^i b_i(S)}  \qquad
\end{eqnarray}
where $b_i(\cdot)$ denotes the $i$-th Betti number, i.e., the rank
of the $i$-th ordinary (co)homology with $\C$-coefficients. 
For a general smooth quasi-projective surface $S$, it is unclear 
whether the above formula still holds. In the following, we let 
$S = \C^2$ and show that the above formula holds for $S = \C^2$.

\begin{proposition}  \label{c2_bi}
$\displaystyle{\sum_{n=0}^{+\infty} \left (\sum_{i} (-1)^i
      b_i\big ( \hil{n, n+1} \big ) z^i\right )q^n
= \frac{1}{1-z^2q} \cdot  \prod_{n=1}^{+\infty}
  \frac{1}{1-z^{2n-2}q^n}}$.
\end{proposition}
\begin{proof}
Let $O$ be the origin of $S = \C^2$. The ${\mathbb
T}_+$-action on $\hil{n, n+1}$ gives rise to a cell decomposition
$\mathcal C_+^n$ of the punctual incidence Hilbert scheme:
\begin{eqnarray*}
S_O^{[n, n+1]} = \big \{ (\xi, \xi') \in \hil{n, n+1}| \Supp(\xi') =
\{O\} \big \}.
\end{eqnarray*}
By the Proposition~2.6.4 in \cite{Ch1}, the dimension of the 
positive part of the tangent space of $\hil{n, n+1}$ at a 
${\mathbb T}_+$-fixed point $\xi_{\la, \mu} = (\xi_\la, \xi_\mu)$ 
is equal to
\begin{eqnarray}  \label{c2_bi.1}
(n+1) - \mu_1,
\end{eqnarray}
where $\mu = (\mu_1 \ge \mu_2 \ge \ldots \ge \mu_r)$. By the 
Theorem~3.3.3~(5) of \cite{Ch1},
\begin{eqnarray}  \label{c2_punc_bi}
\sum_{n=0}^{+\infty} \left (\sum_{i} (-1)^i b_i^{lf} \big ( S_O^{[n,
n+1]} \big ) z^i\right )q^n = \frac{1}{1-z^2q} \cdot
\prod_{n=1}^{+\infty}
  \frac{1}{1-z^{2n-2}q^n}
\end{eqnarray}
where $b_i^{lf}( \cdot )$ stands for the rank of the Borel-Moore
homology group $H_i^{lf}(\cdot)$. On the other hand, the ${\mathbb
T}_-$-action on $\hil{n, n+1}$ gives rise to a cell decomposition
$\mathcal C_-^n$ of $\hil{n, n+1}$ itself. The positive part of the
tangent space of $\hil{n, n+1}$ at a ${\mathbb T}_-$-fixed point
$\xi_{\la, \mu}$ is precisely the negative part of the tangent space
of $\hil{n, n+1}$ at the same ${\mathbb T}_+$-fixed point 
$\xi_{\la, \mu}$. Since $\hil{n, n+1}$ is of dimension $2(n+1)$, 
we see from (\ref{c2_punc_bi}) that
\begin{eqnarray*}
\sum_{n=0}^{+\infty} \left (\sum_{i} (-1)^i b_{4(n+1)-i}^{lf} \big (
S^{[n, n+1]} \big ) z^i\right )q^n = \frac{1}{1-z^2q} \cdot
\prod_{n=1}^{+\infty}
  \frac{1}{1-z^{2n-2}q^n}.
\end{eqnarray*}
Again, since $\hil{n, n+1}$ is smooth, there are natural
isomorphisms:
\begin{eqnarray}   \label{lf_h_coh}
H^{lf}_{4(n+1)-i}(\hil{n, n+1}) \cong H_{i}(\hil{n, n+1})^* \cong
H^i(\hil{n, n+1})
\end{eqnarray}
where $H_{i}(\hil{n, n+1})^*$ is the dual of $H_{i}(\hil{n, n+1})$.
So we obtain the desired formula.
\end{proof}

\subsection{\bf A ${\mathbb T}$-invariant cell decomposition}
\label{subsect_T_cell} $\,$
\par

In the rest of this section, we let $S = \C^2$. Recall that 
the {\it conjugate} of a partition $\mu$ is the partition $\mu'$ 
whose Young diagram is the transpose of that of $\mu$. Then
\begin{eqnarray}   \label{la1_leng}
\ell(\mu) = \mu_1'
\end{eqnarray}
if $\mu' = (\mu_1' \ge \mu_2' \ge \ldots \ge \mu_r')$. In addition,
$\mu$ and $\mu'$ have the same step length.

\begin{proposition}  \label{cell_dec}
Let $S = \C^2$. Then $\hil{n, n+1}$ admits a cell decomposition
\begin{eqnarray}   \label{cell_dec.1}
\hil{n, n+1} \,\, = \,\, \coprod_{\mu \vdash (n+1), \, m_i(\mu') >
0} \,\, S_{\mu, i}
\end{eqnarray}
such that $(\xi_{(\mu')^{(i)}}, \xi_{\mu'}) \in S_{\mu, i} \cong 
\C^{(n+1)+ \ell(\mu)}$
and every cell $S_{\mu, i}$ is ${\mathbb T}$-invariant.
\end{proposition}
\begin{proof}
 From the proof of Proposition~\ref{c2_bi}, we see that
the ${\mathbb T}_-$-action on $\hil{n, n+1}$ gives rise to a cell
decomposition $\mathcal C_-^n$ of $\hil{n, n+1}$. Let $\mathcal
C_-^{\mu', i}$ be the cell corresponding to the ${\mathbb
T}_-$-fixed point $(\xi_{(\mu')^{(i)}}, \xi_{\mu'})$ where 
$m_i(\mu') > 0$. By (\ref{c2_bi.1}) and (\ref{la1_leng}),
\begin{eqnarray*}
\dim (\mathcal C_-^{\mu', i}) = 2(n+1) - [(n+1) - \mu_1'] = (n+1) +
\ell(\mu).
\end{eqnarray*}
Define $S_{\mu, i} = \mathcal C_-^{\mu', i}$ for every partition
$\mu \vdash (n+1)$ and $m_i(\mu') > 0$. Then we have the cell
decomposition (\ref{cell_dec.1}) with $(\xi_{(\mu')^{(i)}}, 
\xi_{\mu'}) \in S_{\mu, i} \cong \C^{(n+1)+\ell(\mu)}$.

To show that $S_{\mu, i} = \mathcal C_-^{\mu', i}$ is ${\mathbb
T}$-invariant, let $(\eta, \eta') \in \mathcal C_-^{\mu', i}$. Then,
\begin{eqnarray*}
\lim_{b \to 0} \, b(\eta, \eta') = (\xi_{(\mu')^{(i)}}, \xi_{\mu'})
\end{eqnarray*}
where $b \in {\mathbb T}_-$. Let $a \in {\mathbb T}$. To show
$a(\eta, \eta') \in \mathcal C_-^{\mu', i}$, it suffices  to verify
\begin{eqnarray}   \label{cell_dec.2}
\lim_{b \to 0} \, b\big (a(\zeta) \big ) = \lim_{b \to 0} \,
b(\zeta)
\end{eqnarray}
for every $\zeta \in S^{[n]}$. Let $f = f(w, z) \in I_\zeta \subset
\C[w,z]$. Then the contribution of $f \in I_\zeta$ to the limiting
ideal $\lim_{b \to 0} I_{b(\zeta)}$ is equal to $w^{i(f)}z^{j(f)}$
where
\begin{eqnarray*}
j(f) &=& \text{max}\{ j|\, \text{ for some $i$,
            $w^iz^j$ is a term in } f\},   \\
i(f) &=& \text{max}\{ i|\, \text{$w^iz^{j(f)}$
            is a term in } f\}.
\end{eqnarray*}
Since $j(a(f)) = j(f)$ and $i(a(f)) = i(f)$, we conclude that the
contribution of $a(f) \in I_{a(\zeta)}$ to $\lim_{b \to 0}
I_{b(a(\zeta))}$ is also $w^{i(f)}z^{j(f)}$. This proves
(\ref{cell_dec.2}).
\end{proof}

\begin{corollary}   \label{cor_e_c}
{\rm (i)} $H_{\mathbb T}^{2k}(S^{[n, n+1]}) = t^{k-(n+1)} \cup
H_{\mathbb T}^{2(n+1)}(S^{[n, n+1]})$ when $k \ge (n+1)$;

{\rm (ii)} There exists a ring isomorphism $H_{\mathbb T}^*(S^{[n,
n+1]})/(t) \cong H^*(S^{[n, n+1]})$.
\end{corollary}
\begin{proof}
(i) Let $\overline{S_{\mu, i}}$ be the closure of $S_{\mu, i}$ in
$\hil{n, n+1}$. By Proposition~\ref{cell_dec},
\begin{eqnarray}       \label{cor_e_c.1}
H_{\mathbb T}^r(S^{[n, n+1]}) \,\, = \,\, \bigoplus_{\mu \vdash
(n+1), \, m_i(\mu') > 0 \atop 2j + 2(n+1) - 2\ell(\mu) = r} \,\, \C
t^j \cup  \big [\, \overline{S_{\mu, i}} \, \big ]
\end{eqnarray}
Here and below, $[\cdot]$ denotes the equivariant fundamental cycle
or its associated equivariant cohomology class. Now (i) follows
immediately.

(ii) There is the forgetful map $H_{\mathbb T}^*(S^{[n, n+1]}) \to
H^*(S^{[n, n+1]})$ which is a ring homomorphism. By
Proposition~\ref{cell_dec} and (\ref{lf_h_coh}), a $\C$-linear basis
of $H^*(S^{[n, n+1]})$ consists of the (ordinary) fundamental
cohomology classes of the closures $\overline{S_{\mu, i}}$.
Combining this with (\ref{cor_e_c.1}), we obtain $H_{\mathbb
T}^*(S^{[n, n+1]})/(t) \cong H^*(S^{[n, n+1]})$.
\end{proof}
\subsection{The equivariant Zariski tangent spaces of $\hil{n}$}
\label{subsect_zts_hil} $\,$
\par
The first study of the equivariant Zariski tangent space of the
Hilbert scheme $\hil{n}$ at the fixed points was carried out in 
\cite{ES1}. Here we review an approach due to Cheah in \cite{Ch1}. 
It will be used in the next subsection for the equivariant 
Zariski tangent space of the incidence Hilbert scheme $\hil{n,n+1}$.

Let $\la \vdash n$. Then the ${\mathbb T}$-invariant ideal
$I_{\xi_\la} \subset R = \C[w, z]$ is given by:
\begin{eqnarray*}
I_{\xi_{\la}} = (w^{\la_1}, zw^{\la_2}, \ldots, z^{r-1}w^{\la_r},
z^r).
\end{eqnarray*}
The Zariski tangent space $T_{\xi_{\la}}S^{[n]}$ of $S^{[n]}$ at
$\xi_{\la}$ is canonically isomorphic to the space
$\Hom(I_{\xi_{\la}}, R/I_{\xi_{\la}})$. To obtain a pure-weight
linear basis of the ${\mathbb T}$-invariant space
$\Hom(I_{\xi_{\la}}, R/I_{\xi_{\la}})$, we represent $I_{\xi_{\la}}$
by a Young diagram $D_\la$ as in \cite{Ch1}.

\begin{example}
The ideal $I = \langle w^5, zw^4, z^4w^2, z^6 \rangle$ is
represented by the diagram:
\begin{figure}
[H]{ \centerline{\includegraphics{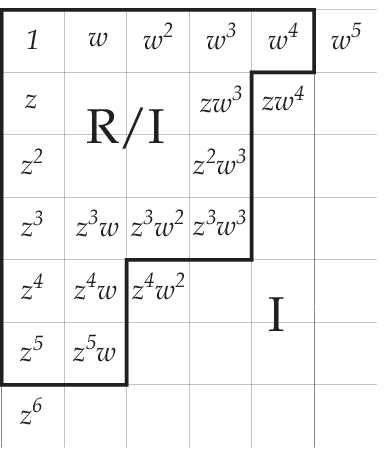}}}\caption{ }
\end{figure}
\end{example}

If we look at the Young diagram of $I = I_{\xi_\la}$, then the
corners of its complement (that is, the shaded boxes in the second
diagram of Figure 2) represent the unique minimal set of monomials
that generate the ideal $I_{\xi_\la}$. Denote this set by $A$, and
let $B$ be the set of monomials not in $I_{\xi_\la}$:
\begin{figure}
[H]{\centerline{\includegraphics{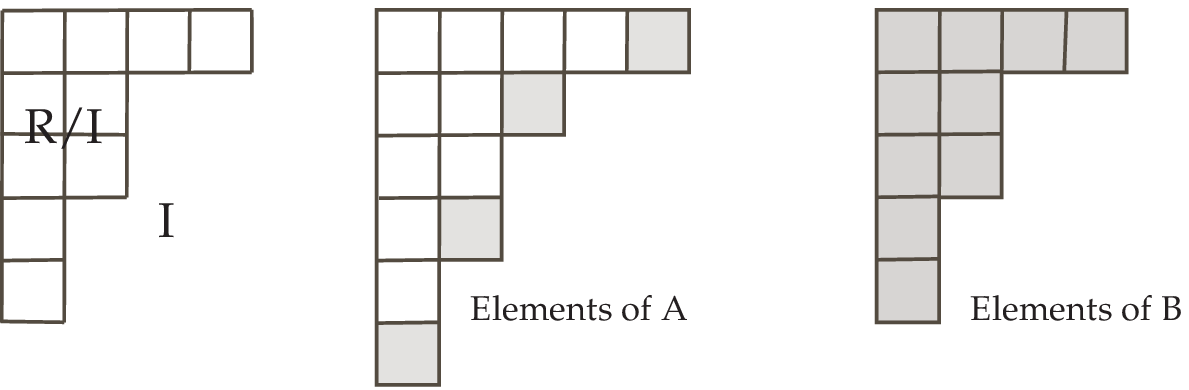}}}\caption{ }
\end{figure}

Rewrite the partition $\la = (\la_1 \ge \la_2 \ge \ldots \ge \la_r)$
as
\begin{eqnarray}   \label{p_i}
& &\la_1 = \ldots = \la_{p_0}   \nonumber\\
&>&\la_{p_0+1} = \ldots = \la_{p_0+p_1}  \nonumber \\
&>&\ldots  \,\, >           \nonumber  \\
&>&\la_{p_0+p_1+\ldots+p_{m-2}+1} = \ldots
   = \la_{p_0+p_1+\ldots+p_{m-2}+p_{m-1}} = \la_r.
\end{eqnarray}
Recall the step length $s(\la)$ from Definition~\ref{step}~(i). We
have
\begin{eqnarray}   \label{ms(la)}
m = s(\la).
\end{eqnarray}
The elements (called the {\it canonical generators}) in the set $A$
are:
\begin{eqnarray} \label{alpha_i}
\alpha_0 &:=& w^{\la_{p_0}}, \,\,   \nonumber   \\
\alpha_1 &:=& z^{p_0}w^{\la_{p_0+p_1}}, \nonumber  \\
&\ldots,&  \nonumber  \\
\alpha_{m-1} &:=& z^{p_0+p_1+\ldots+p_{m-2}}
  w^{\la_{p_0+p_1+ \ldots+p_{m-2}+p_{m-1}}}, \nonumber \\
\alpha_{m} &:=& z^{p_0+p_1+\ldots+p_{m-2}+p_{m-1}} = z^{r}.
\end{eqnarray}
Define $q_m = \la_{p_0+p_1+ \ldots+p_{m-2}+p_{m-1}}$. For $1 \le i
\le m-1$, define
\begin{eqnarray}     \label{q_i}
q_i = \la_{p_0+p_1+ \ldots +p_{i-1}} - \la_{p_0+p_1+ \ldots
+p_{i-1}+p_{i}}.
\end{eqnarray}
Note that $p_i$ is the vertical distance between the cells
representing $\alpha_i$ and $\alpha_{i+1}$ and that $q_i$ is 
the horizontal distance between the cells representing 
$\alpha_i$ and $\alpha_{i-1}$.

 For $\alpha = \alpha_i \in A$, let $P_\alpha$ be the
subset of $B$ consisting of the elements $b$ satisfying
\par
$\quad$ (i) $b$ lies to the left of $\alpha$ in the Young diagram,
\par
$\quad$ (ii) $z^{p_i}b \in I_{\xi_\la}$,
\par\noindent
and let $Q_\alpha$ be the subset of $B$ consisting of the elements
$b$ satisfying
\par
$\quad$ (i) $b$ lies above $\alpha$ in the Young diagram,
\par
$\quad$ (ii) $w^{q_i}b \in I_{\xi_\la}$.

Let $\mathcal S$ be the subset of $\Hom(I_{\xi_{\la}},
R/I_{\xi_{\la}})$ consisting of elements of pure weight which take
canonical generators in $A$ either to zero or to monomials in $B$
modulo $I_{\xi_\la}$. For $\beta \in P_\alpha \cup Q_\alpha$, define
$f_{\alpha, \beta} \in \mathcal S$ to be the unique element
satisfying
\par
$\quad$ (i) $f_{\alpha, \beta}(\alpha) = \beta$,
\par
$\quad$ (ii) $f_{\alpha, \beta}$ takes the largest number
             of canonical generators to zero.
\par\noindent
The conditions (i) and (ii) imply that $f_{\alpha_i,
\beta}(\alpha_j) = 0$ if $\beta \in P_{\alpha_i}$ and $j > i$, and
that $f_{\alpha_i, \beta} (\alpha_j) = 0$ if $\beta \in
Q_{\alpha_i}$ and $j < i$. By the Proposition~2.5.4 of \cite{Ch1}, a
pure weight basis of the tangent space $T_{\xi_{\la}}S^{[n]}
\cong\Hom(I_{\xi_{\la}}, R/I_{\xi_{\la}})$ of $\hil{n}$ at $\xi_\la$
is
\begin{eqnarray*}
\{f_{\alpha, \beta}|\,\, \alpha \in A, \, \beta \in P_\alpha \cup
Q_\alpha \}.
\end{eqnarray*}
When $\beta \in Q_\alpha$, we have $\beta = \alpha \cdot
w^{i_1}/z^{i_2}$ for some integers $i_1 \ge 0$ and $i_2 > 0$, and
the weight of $f_{\alpha, \beta}$ is equal to $(i_1+i_2)$, which is
the hook length $h(\square)$ of certain cell $\square$ in the Young
diagram $D_\la$. As $\beta$ runs in the set $Q_\alpha$, $\square$
runs over all the cells in $D_\la$ exactly once. Similarly, when
$\beta \in P_\alpha$, we have $\beta = \alpha \cdot z^{i_2}/w^{i_1}$
for some integers $i_1 > 0$ and $i_2 \ge 0$, and the weight of
$f_{\alpha, \beta}$ is equal to $-(i_1+i_2)$, where $(i_1+i_2)$ is
the hook length $h(\square)$ of certain cell $\square$ in $D_\la$.
Again, as $\beta$ runs in $Q_\alpha$, $\square$ runs over all the
cells in $D_\la$ exactly once. Hence, there exists a $\mathbb
T$-equivariant identification:
\begin{eqnarray}   \label{hil_ts.1}
T_{\xi_{\la}}S^{[n]} = \bigoplus_{\square \in D_\la} \big (
\theta^{h(\square)} \oplus \theta^{-h(\square)} \big ).
\end{eqnarray}
It follows that the $\mathbb T$-equivariant Euler class of the
tangent space is
\begin{eqnarray}   \label{hil_ts.2}
e_{\mathbb T} \left (T_{\xi_{\la}}S^{[n]} \right ) = (-1)^n \cdot
\prod_{\square \in D_\la} h(\square)^2 \cdot t^{2n} = (-1)^n \cdot
h(\la)^2 \cdot t^{2n}
\end{eqnarray}
where $h(\la)$ is the product of all the hook lengths $h(\square),
\square \in D_\la$.
\subsection{The equivariant Zariski tangent spaces of $\hil{n, n+1}$}
\label{subsect_zts_inc_hil} $\,$
\par

Let $(\la, \mu)$ be an incidence pair of partitions with $\la \vdash
n$. Then, $\big (\xi_{\la}, \xi_{\mu} \big)$ is a $\mathbb T$-fixed
point in $\hil{n, n+1}$. There are $\mathbb T$-equivariant maps:
\begin{eqnarray}
&\phi:&\Hom \big (I_{\xi_{\la}}, R/I_{\xi_{\la}}
  \big ) \to \Hom \big (I_{\xi_{\mu}}, R/I_{\xi_{\la}}
  \big ),  \label{tang.1}  \\
&\psi:&\Hom \big (I_{\xi_{\mu}}, R/I_{\xi_{\mu}} \big )
  \to \Hom \big (I_{\xi_{\mu}}, R/I_{\xi_{\la}} \big ).
  \label{tang.2}
\end{eqnarray}
 From pages 42-43 in \cite{Ch1}, we see that the
Zariski tangent space of $\hil{n, n+1}$ at the point $\big
(\xi_{\la}, \xi_{\mu} \big)$ is canonically isomorphic to
$\ker(\phi-\psi)$ where
\begin{eqnarray*}
(\phi-\psi):\,\, \Hom \big (I_{\xi_{\la}}, R/I_{\xi_{\la}} \big )
\oplus \Hom \big (I_{\xi_{\mu}}, R/I_{\xi_{\mu}} \big ) \to \Hom
\big (I_{\xi_{\mu}}, R/I_{\xi_{\la}} \big )
\end{eqnarray*}
is defined by letting $(\phi-\psi)(a, b) = \phi(a) - \psi(b)$. By
the Lemma~2.6.2 in \cite{Ch1}, $(\phi-\psi)$ is surjective.
Therefore, there is a $\mathbb T$-equivariant exact sequence:
\begin{eqnarray}   \label{tang.3}
0 \to \ker(\phi-\psi) \to \Hom \big (I_{\xi_{\la}}, R/I_{\xi_{\la}}
\big ) \oplus \Hom \big (I_{\xi_{\mu}},
R/I_{\xi_{\mu}} \big ) \qquad \nonumber  \\
\to \Hom \big (I_{\xi_{\mu}}, R/I_{\xi_{\la}} \big ) \to 0.
\end{eqnarray}

We keep using the notations $A, B, p_i, \alpha_i, q_i$ associated to
$\la$ from \S \ref{subsect_zts_hil}, and let $A', B'$ be the
corresponding notations associated to $\mu$. Put
\begin{eqnarray*}
A' = \{\alpha_0', \alpha_1', \ldots, \alpha_s'\}.
\end{eqnarray*}
The Young diagram of $I_{\xi_{\mu}}$ is obtained from that
of $I_{\xi_{\la}}$ by adding one of the cells which represents a
canonical generator in $A$. Let
\begin{eqnarray}   \label{alpha_k}
\alpha_k \in A
\end{eqnarray}
be this canonical generator. Then, $\alpha_k \in B'$ and so
$\alpha_k \not \in A'$. Note that
\begin{eqnarray*}
\alpha_k \in P_{\alpha_i'} \cup Q_{\alpha_i'}
\end{eqnarray*}
for all $0 \le i \le s$, and that the homomorphism $f_{\alpha_i',
\alpha_k} \in \Hom \big (I_{\xi_{\mu}}, R/I_{\xi_{\mu}} \big )$ maps
$\alpha_i'$ to $\alpha_k$ and all the other canonical generators in
$A'$ to zero. Moreover, $f_{\alpha_0', \alpha_k}, \ldots,
f_{\alpha_s', \alpha_k}$ form a basis of $\ker(\psi) \subset
\ker(\phi-\psi)$, i.e., we have
\begin{eqnarray}   \label{fai'ak}
\ker(\phi-\psi) \supset \ker(\psi) = \bigoplus_{i=0}^s \C
f_{\alpha_i', \alpha_k}.
\end{eqnarray}

\begin{definition}  \label{squareii'}
Let $(\la, \mu)$ be an incidence pair, and $k$ be from
(\ref{alpha_k}).
\begin{enumerate}
\item[{\rm (i)}] For $0 \le i \le k-1$, let $\square_{k, i}$ be
the cell in the Young diagram $D_\la$ which is directly to the left
of $\alpha_i$ and directly above $\alpha_k$, and let $\square_{k,
i}'$ be the cell which is the $(p_i-1)$-th cell directly under
$\square_{k, i}$ ($\square_{k, i}' = \square_{k, i}$ if $p_i = 1$).
For $k+1 \le i \le m = s(\la)$, let $\square_{k, i}$ be the cell
which is directly above $\alpha_i$ and directly to the left of
$\alpha_k$, and let $\square_{k, i}'$ be the cell which is the
$(q_{i}-1)$-th cell directly to the right of $\square_{k, i}$
($\square_{k, i}' = \square_{k, i}$ if $q_i = 1$).

\item[{\rm (ii)}] Define $k(\la, \mu)$ to be the integer
$k$, and define
\begin{eqnarray}   \label{hlagamma}
h(\la, \mu) = h(\la)^2 \cdot \prod_{0 \le i \le s(\la) \atop i \ne
k(\la, \mu)} \frac{1+h(\square_{k(\la, \mu), i})}{h(\square_{k(\la,
\mu), i}')}.
\end{eqnarray}
\end{enumerate}
\end{definition}

The following lemma is the key step for determining the ring
structure of the equivariant cohomology ring of the incidence
Hilbert scheme $\hil{n, n+1}$. Since the proof is a bit technical,
we place it in the Appendix.

\begin{lemma} \label{lemma_eT}
Let $\big (\xi_{\la}, \xi_{\mu} \big) \in \hil{n, n+1}$ be a
$\mathbb T$-fixed point. Then the $\mathbb T$-equivariant Euler
class of the tangent space of $\hil{n, n+1}$ at $\big (\xi_{\la},
\xi_{\mu} \big)$ is equal to
\begin{eqnarray}   \label{lemma_eT.0}
e_{\mathbb T} =(-1)^{n+1} h(\la, \mu) \cdot t^{2(n+1)}.
\end{eqnarray}
\end{lemma}

\subsection{A bilinear pairing}
\label{subsect_bp} $\,$
\par

Recall that the ${\mathbb T}$-fixed points in $\hil{n, n+1}$ 
are of the form $\xi_{\la, \mu} = (\xi_\la, \xi_\mu)$
where $(\la, \mu)$ denotes incidence pairs of partitions.
Let
\begin{eqnarray}   \label{inc_mu_i}
\iota_{\la, \mu}: \,\, \xi_{\la, \mu} = (\xi_\la, \xi_\mu)
\hookrightarrow S^{[n, n+1]}
\end{eqnarray}
be the inclusion map. Let $1_{\xi_{\la, \mu}} \in H^0_{\mathbb T}
(\xi_{\la, \mu})$ be the unit. Thus,
\begin{eqnarray*}
[\xi_{\la, \mu}]= (\iota_{\la, \mu})_! (1_{\xi_{\la, \mu}}) 
\in H^{4(n+1)}_{\mathbb T}(S^{[n, n+1]}).
\end{eqnarray*}

Denote by $\mathbb C[t]'$ the localization of the ring $\mathbb
C[t]$ at the ideal $(t-1)$, and denote
\begin{eqnarray*}
\w \iota_n = \bigoplus\limits_{(\la, \mu) \,\, \text{\rm incidence}} 
\iota_{\la, \mu}: \,\, (\hil{n, n+1})^{\mathbb T} \to \hil{n, n+1}.
\end{eqnarray*}
We define $H^*_{\mathbb T}\big ((S^{[n, n+1]})^{\mathbb T} \big
)'=H^*_{\mathbb T}\big ((S^{[n, n+1]})^{\mathbb T}\big )
\otimes_{\mathbb C[t]} \mathbb C[t]'$ and define $H^*_{\mathbb
T}(S^{[n, n+1]} )'$ similarly. Then we have the induced Gysin map:
\begin{eqnarray*}
\w \iota_{n!} \colon H^*_{\mathbb T}\big ( (S^{[n, n+1]})^{\mathbb
T}\big )' \longrightarrow H^*_{\mathbb T}(S^{[n, n+1]} )'.
\end{eqnarray*}
By the localization theorem, $\w \iota_{n!}$ is an isomorphism. The
inverse $(\w \iota_{n!})^{-1}$ is given by
\begin{eqnarray*}
\W A \,\, \mapsto \,\, \left ( \frac{(\iota_{\la, \mu})^*(\W A)}
{e_{\mathbb T}\big (T_{\xi_{\la, \mu}}\hil{n,n+1}\big )} \right
)_{(\la, \mu) \,\, \text{\rm incidence}}.
\end{eqnarray*}
Therefore, we conclude from Lemma~\ref{lemma_eT} on the Euler
class $e_{\mathbb T}$ that
\begin{eqnarray}   \label{loc_for}
(\w \iota_{n!})^{-1}(\W A) \,\, = \,\, \left ( \frac{(\iota_{\la, 
\mu})^*(\W A)}{(-1)^{n+1} h(\la, \mu) \, t^{2(n+1)}} \right
)_{(\la, \mu) \,\, \text{\rm incidence}}.
\end{eqnarray}

Next, we define a bilinear pairing on $H^*_{\mathbb T}(S^{[n,
n+1]})'$ by:
\begin{eqnarray}
&\langle -, - \rangle:
   H^*_{\mathbb T}(S^{[n, n+1]})' \otimes_{\mathbb C[t]'}
   H^*_{\mathbb T}(S^{[n, n+1]})'\to \mathbb C[t]', &
   \label{pairing.1} \\
&\langle \W A, \W B\rangle = (-1)^{n+1} \,
   \w \pi_{n!}(\w \iota_{n!})^{-1}(\W A \cup \W B)&
   \label{pairing.2}
\end{eqnarray}
where $\w \pi_n$ is the projection of the set $(S^{[n,
n+1]})^{\mathbb T}$ of ${\mathbb T}$-fixed points to a point. This
induces a bilinear pairing, again denoted by $\langle -, - \rangle$,
on the space:
\begin{eqnarray}      \label{Wft'}
\Wft' = \bigoplus_{n=0}^{+\infty} H^*_{\mathbb T}(S^{[n, n+1]})'.
\end{eqnarray}

\subsection{A new ring $\Wh_{\mathbb T, n}$ and its linear basis 
from the fixed points}
\label{subsect_wht} 
$\,$
\par

For $n \ge 0$, let $\Wh_{\mathbb T, n} = H_{\mathbb T}^{2(n+1)}
(S^{[n, n+1]})$ be the middle-degree equivariant cohomology of
the incidence Hilbert scheme $\hil{n, n+1}$.  
By Corollary~\ref{cor_e_c}~(i),
\begin{eqnarray*}
H_{\mathbb T}^{4(n+1)}(S^{[n, n+1]}) = t^{(n+1)} \cup H_{\mathbb
T}^{2(n+1)}(S^{[n, n+1]}).
\end{eqnarray*}
Note that $H_{\mathbb T}^*(S^{[n, n+1]})$ is $\C[t]$-torsion free.
Define a product $\w \star$ on $\Wh_{\mathbb T, n}$ by:
\begin{eqnarray}      \label{w*}
t^{n+1} \cup(\W A \w \star \W B)=\W A \cup \W B
\end{eqnarray}
for $\W A, \W B \in \Wh_{\mathbb T, n}$. Then, we see that $\big (
\Wh_{\mathbb T, n}, \w \star \big )$ is a ring.

Next, we construct a linear basis of $\Wh_{\mathbb T, n}$. Define
the class $[{\la, \mu}] \in \Wh_{\mathbb T, n}$ by
\begin{eqnarray}      \label{normalize}
t^{n+1} \cup [{\la, \mu}] = [\xi_{\la, \mu}]
\end{eqnarray}
since $[\xi_{\la, \mu}] \in H^{4(n+1)}_{\mathbb T}(S^{[n, n+1]})$.
Note that we have
\begin{eqnarray} \label{cup-prod}
   [\xi_{\la, \mu}] \cup [\xi_{\w \la, \w \mu}]
&=&(\iota_{\la, \mu})_! (1_{\xi_{\la, \mu}}) \cup
   (\iota_{\w \la, \w \mu})_! (1_{\xi_{\w \la, \w \mu}}) 
   \nonumber \\
&=&(\iota_{\la, \mu})_! \big (1_{\xi_{\la, \mu}}
   \cup \iota_{\la, \mu}^*(\iota_{\w \la, \w \mu})_!
   (1_{\xi_{\w \la, \w \mu}}) \big )  \nonumber \\
&=&\delta_{(\la, \mu),(\w \la, \w \mu)} \,
   e_{\mathbb T}(T_{\xi_{\la, \mu}}\hil{n, n+1})
   [\xi_{\la, \mu}]     \nonumber \\
&=&\delta_{(\la, \mu),(\w \la, \w \mu)} \, (-1)^{n+1} \,
   h(\la, \mu)t^{2(n+1)} [\xi_{\la, \mu}]
\end{eqnarray}
by the projection formula and Lemma \ref{lemma_eT}. Thus we obtain
\begin{eqnarray} \label{w*-prod}
  [\la, \mu] \, \w \star \, [\w \la, \w \mu]
= \delta_{(\la, \mu),(\w \la, \w \mu)} \, (-1)^{n+1} \,
   h(\la, \mu) \,\, [\la, \mu].
\end{eqnarray}
Combining this with the localization theorem, we see that the
classes
\begin{eqnarray*}
[\la, \mu]
\end{eqnarray*}
where $\mu \vdash (n+1)$ and $(\la, \mu)$ is an incidence pair,
form a linear basis of $\Wh_{\mathbb T, n}$.

In addition, we obtain from (\ref{pairing.2}) and (\ref{cup-prod})
that
\begin{eqnarray}      \label{muinuj}
   \big \langle [{\la, \mu}], [\w \la, \w \mu] \big \rangle
&=&(-1)^{n+1} \, \w \pi_{n!}(\w \iota_{n!})^{-1}\big ([\la, \mu]
     \cup [\w \la, \w \mu] \big )    \nonumber   \\
&=&(-1)^{n+1} \, \w \pi_{n!}(\w \iota_{n!})^{-1}\big (t^{-2(n+1)}
     \cup [\xi_{\la, \mu}] \cup[\xi_{\w \la, \w \mu}] \big )
     \nonumber    \\
&=&\delta_{(\la, \mu),(\w \la, \w \mu)} \,
     h(\la, \mu) \cdot \w \pi_{n!}(\w \iota_{n!})^{-1}
     \big ([\xi_{\la, \mu}]\big ) \nonumber   \\
&=&\delta_{(\la, \mu),(\w \la, \w \mu)} \, h(\la, \mu).
\end{eqnarray}
It follows that the restriction to $\Wh_{\mathbb T, n}$ of the
bilinear form $\langle -, -\rangle$ on the space $H^*_{\mathbb
T}(S^{[n, n+1]})'$ is a nondegenerate bilinear form:
\begin{eqnarray} \label{bil_wh}
\langle -, - \rangle:\,\, \Wh_{\mathbb T, n} \times \Wh_{\mathbb T,
n} \to \C
\end{eqnarray}
This induces a nondegenerate bilinear form, denoted again by
$\langle -, -\rangle$, on:
\begin{eqnarray}      \label{Wfte}
\Wft = \bigoplus_{n=0}^{+\infty} \Wh_{\mathbb T, n}.
\end{eqnarray}
Note that a $\mathbb C$-linear basis of the vector space
$\Wh_{\mathbb T}$ is given by
\begin{eqnarray}   \label{WB1}
  \W {\mathcal B}_1
= \big \{ [\la, \mu] \big \}_{(\la, \mu) \,\, \text{\rm incidence}}.
\end{eqnarray}

\section{\bf The loop algebra action on $\Wft$}
\label{sect_loop}

One of the most important features of the Hilbert schemes $\hil{n}$ 
of points on a surface $S$ is a Heisenberg algebra action on 
the direct sum of the cohomology groups of $\hil{n}$ over all $n$ 
discovered by Nakajima and Grojnowski \cite{Gro, Na1}. 
It lays the foundation for a new method in the study of 
the cohomology ring of the Hilbert scheme $\hil{n}$. 
Without much difficulty, one can transport the Heisenberg algebra 
action to the equivariant cohomology of $\hil{n}$ when $S = \C^2$.

A loop algebra of a Heisenberg algebra was found in \cite{LQ} to act
on the direct sum of the cohomology groups of incidence Hilbert
schemes $\hil{n,n +1}$. 
It is generated by a Heisenberg algebra and a translation operator. 
In this section, we transport the results in \cite{LQ} to the 
equivariant cohomology of $\hil{n,n+1}$ when $S = \C^2$.

\subsection{The Heisenberg operators}
\label{subsect_hei} $\,$
\par

In the rest of this section, let $S = \C^2$.
Let $C^w$ and $C^z$ be the $w$-axis and $z$-axis of $S = \C^2$
respectively. By the localization theorem, we have
\begin{eqnarray} \label{Cwz}
[C^w] = t = -t^{-1}[O], \qquad [C^z] = -t = t^{-1}[O]
\end{eqnarray}
in $H^2_{\mathbb T}(S)$, where $O \in S$ is the origin. In
particular, $[C^w] = -[C^z]$ in $H^*_{\mathbb T}(S)$.

Next, let $Y = C^w$ or $C^z$. Then, $Y$ is $\mathbb T$-invariant.
For $m \ge 0$ and $n > 0$, define the closed subset
$\Wq^{[m+n,m]}_Y$ of $\hil{m+n, m+n+1} \times \hil{m, m+1}$:
\begin{eqnarray*}
\Wq^{[m+n,m]}_Y= \big \{ \big ( (\xi,\xi'), (\eta,\eta') \big ) |
&&\xi\supset \eta, \,\, \xi' \supset \eta', \,\,
  \Supp(I_\eta/I_\xi)= \{ s \} \subset Y, \\
&&\Supp(I_\xi/I_{\xi'}) = \Supp(I_\eta/I_{\eta'}) \big \}.
\end{eqnarray*}
Then $\Wq^{[m+n,m]}_Y$ is $\mathbb T$-invariant. Define the linear
operator $\wa_{-n}([Y]) \in \End(\Wft')$ by
\begin{eqnarray}   \label{waY}
\wa_{-n}([Y])(\W A) = D^{-1} \w p_{1!} \left ( \w p_2^*\W A \cap
\big [\Wq^{[m+n,m]}_Y \big ]\right )
\end{eqnarray}
for $\W A \in H^*_{\mathbb T}(S^{[m, m+1]})'$, where $\w p_1$ and
$\w p_2$ are the two projections of $\hil{m+n, m+n+1} \times \hil{m,
m+1}$. Note that the restriction of $\w p_1$ to $\Wq^{[m+n,m]}_Y$ is
proper. Define $\wa_{n}([Y]) \in \End(\Wft')$ to be the adjoint
operator of $\wa_{-n}([Y])$ with respect to the bilinear form
$\langle -, - \rangle$ on the linear space $\Wft'$. Alternatively, 
we have
\begin{eqnarray}
   \wa_{-n}([Y])(\W A)
&=&D^{-1} (\w p_1')_! (\text{Id}_{\hil{m+n, m+n+1}} \times
     \w \iota_{m})_!^{-1} \left ( \w p_2^*\W A \cap
     \big [\Wq^{[m+n,m]}_Y \big ]\right ),   \label{alt1}   \\
   \wa_{n}([Y])(\W B)
&=&(-1)^n \, D^{-1}(\w p_2')_! (\w \iota_{m+n} \times
     \text{Id}_{\hil{m, m+1}})_!^{-1} \left (\w p_1^*\W B \cap
     \big [\Wq^{[m+n,m]}_Y \big ]\right ) \qquad    \label{alt2}
\end{eqnarray}
for $\W A \in H^*_{\mathbb T}(S^{[m, m+1]})'$ and $\W B \in
H^*_{\mathbb T}(S^{[m+n, m+n+1]})'$, where $\w p_1', \w p_2'$ are
the projections:
\begin{eqnarray*}
&\w p_1':&\hil{m+n, m+n+1} \times \big ( \hil{m, m+1}
          \big )^{\mathbb T} \to \hil{m+n, m+n+1}, \\
&\w p_2':&\big ( \hil{m+n, m+n+1} \big )^{\mathbb T} \times
          \hil{m, m+1} \to \hil{m, m+1}.
\end{eqnarray*}

Let $n > 0$. From the definition of $\wa_{-n}([Y]) \in \End(\Wft')$,
we see that
\begin{eqnarray*}
\wa_{-n}([Y])(\W A) \in \Wh_{\mathbb T, m+n}
\end{eqnarray*}
if $\W A \in \Wh_{\mathbb T, m} = H^{2(m+1)}_{\mathbb T}(S^{[m,
m+1]}) \subset H^*_{\mathbb T}(S^{[m, m+1]})'$. Hence the
restriction of $\wa_{-n}([Y])$ to the space $\Wft$ gives 
a linear operator in $\End(\Wft)$, denoted by $\wa_{-n}([Y])$ 
as well. Recall from (\ref{bil_wh}) that there is a bilinear form
\begin{eqnarray*}
\langle -, -\rangle \colon \Wft \otimes_{\mathbb C} \Wft \to \mathbb
C,
\end{eqnarray*}
which is the restriction of the bilinear form $\langle -, -\rangle$
on $\Wft'$. Thus, the restriction of $\wa_n([Y])$ to $\Wft$ 
is the adjoint operator of $\wa_{-n}([Y])$ with respect to the bilinear
form $\langle-, -\rangle$ on $\Wft$, and hence is an operator in
$\End(\Wft)$ which will again be denoted by $\wa_n([Y])$. Finally,
we define $\wa_0([Y]) = 0 \in \End(\Wft)$.

\begin{proposition}   \label{prop_hei_comm}
The operators $\wa_n([C^z])$, $n \in \mathbb Z$, acting on the space
$\Wft$ satisfy the following Heisenberg commutation relation:
\begin{eqnarray}        \label{Heisenberg}
\big [\wa_n([C^z]), \wa_m([C^z]) \big ] = n \delta_{n, -m} \,\,
\text{\rm Id}_{\Wft}.
\end{eqnarray}
\end{proposition}
\begin{proof}
Since $[C^w] = -[C^z]$, we have $[\Wq^{[m+n,m]}_{C^w}] =
-[\Wq^{[m+n,m]}_{C^z}]$. It follows from the definition that
$\wa_m([C^w]) = -\wa_m([C^z])$. Hence (\ref{Heisenberg}) is
equivalent to
\begin{eqnarray}        \label{Heisenberg.1}
\big [\wa_n([C^w]), \wa_m([C^z]) \big ] = -n \delta_{n, -m} \,\,
\text{\rm Id}_{\Wft}.
\end{eqnarray}
Note that $C^w$ and $C^z$ intersect transversely at the origin. By
(\ref{alt1}) and (\ref{alt2}), the commutation relation
(\ref{Heisenberg.1}) is reduced to the intersections between certain
cycles related to various $\big [\Wq^{[\ell_1, \ell_2]}_{C^w} \big
]$ and certain cycles related to various $\big [\Wq^{[\ell_3,
\ell_4]}_{C^z} \big ]$. Therefore, an argument similar to the one
used in the proof of the Proposition~3.5 in \cite{LQ} (for smooth
projective surfaces) works in our situation. This proves
(\ref{Heisenberg.1}).
\end{proof}

\subsection{The translation operator}
\label{subsect_trans} $\,$
\par

For $m \ge 0$, define $\Wq_m \subset S^{[m+1, m+2]} \times
S^{[m,m+1]}$ to be the closed subset:
\begin{eqnarray*}
  \Wq_m
= \{\big ( (\xi^{\prime}, \xi^{\prime \prime}), (\xi,
   \xi^{\prime}) \big ) |\Supp(I_\xi/I_{\xi'}) =
   \Supp(I_{\xi'}/I_{\xi''})\}.
\end{eqnarray*}
Then, the subset $\Wq_m$ is $\mathbb T$-invariant, and $\dim \Wq_m =
(2m + 3)$.

\begin{definition}  \label{def_oper_wt}
Define the linear operator $\wt^{\mathbb T} \in \End \big (\Wft \big
)$ by
\begin{eqnarray}  \label{def_oper_wt.1}
\wt^{\mathbb T}(\W A) = D^{-1}\tilde{p}_{1*}(\tilde{p}_2^*\W A \cap
[\Wq_m])
\end{eqnarray}
for $\W A \in \Wh_{\mathbb T, m}$, where $\tilde{p}_1, \tilde{p}_2$
are the two projections of $S^{[m+1, m+2]} \times S^{[m,m+1]}$.
\end{definition}

\begin{proposition}  \label{comm_wtwa}
{\rm (i)} The adjoint operator $\big (\wt^{\mathbb T} \big
)^\dagger$ is the left inverse of $\wt^{\mathbb T}$;

{\rm (ii)} $\wt^{\mathbb T}$ and $\big (\wt^{\mathbb T} \big
)^\dagger$ commute with the Heisenberg operators $\wa_{-n}([C^z])$.
\end{proposition}
\begin{proof}
Note that $\wt$ and its adjoint operator $\wt^\dagger$ are also
given by
\begin{eqnarray*}
   \wt^{\mathbb T}(\W A)
&=&D^{-1} (\w p_1')_! (\text{Id}_{\hil{m+1, m+2}} \times
     \w \iota_{m})_!^{-1} \left ( \w p_2^*\W A \cap [\Wq_m]\right ),\\
   \big (\wt^{\mathbb T} \big )^\dagger(\W B)
&=&-D^{-1}(\w p_2')_! (\w \iota_{m+1} \times
     \text{Id}_{\hil{m, m+1}})_!^{-1} \left (\w p_1^*\W B \cap
     [\Wq_m] \right )
\end{eqnarray*}
for $\W A \in \Wh_{\mathbb T, m}$ and $\W B \in \Wh_{\mathbb T,
m+1}$, where $\w p_1'$ and $\w p_2'$ are the projections:
\begin{eqnarray*}
&\w p_1':&\hil{m+1, m+2} \times \big ( \hil{m, m+1} \big
            )^{\mathbb T} \to \hil{m+1, m+2}, \\
&\w p_2':&\big ( \hil{m+1, m+2} \big )^{\mathbb T} \times
            \hil{m, m+1} \to \hil{m, m+1}.
\end{eqnarray*}
So our results follow from arguments similar to the proofs of the
Lemma~4.2~(i) and Proposition~4.3 in \cite{LQ} for smooth projective
surfaces.
\end{proof}

\subsection{The loop algebra action}
\label{subsect_loop} $\,$
\par
\begin{definition}   \label{alg_h}
Let $\wa_n^{\mathbb T} = \wa_n([C^z])$ for all $n \in \Z$. Define
${\w {\mathfrak h}}_{\mathbb T}$ to be the Heisenberg algebra
generated by the operators $\wa_n^{\mathbb T}$, $n \in \Z$, and the
identity operator ${\rm Id}_{\Wft}$.
\end{definition}

The loop algebra of a Lie algebra $\mathfrak g$ is $\C [ u, u^{-1}]
\otimes_\C \mathfrak g$ with
\begin{eqnarray}   \label{loop}
[u^m \otimes g_1,  u^n \otimes g_2] =  u^{m+n} \otimes [g_1, g_2].
\end{eqnarray}

\begin{theorem}  \label{thm_struc}
The space $\Wft$ is a representation of the Lie algebra
$\C[u^{-1}]\otimes_\C \w{\mathfrak h}_{\mathbb T}$ with a highest
weight vector being the vacuum vector
\begin{eqnarray}   \label{thm_struc.1}
\vac = [C^z] \in \Wh_{\mathbb T, 0} = H_{\mathbb T}^2(\hil{0,1}) =
H_{\mathbb T}^2(S) = H_{\mathbb T}^2(\C^2)
\end{eqnarray}
where $u^{-1}$ acts via $\wt^{\mathbb T}$, and $C^z$ denotes the
$z$-axis of $S = \C^2$.
\end{theorem}
\begin{proof}
Follows from Proposition~\ref{c2_bi},
Proposition~\ref{prop_hei_comm} and Proposition~\ref{comm_wtwa}.
\end{proof}

For a partition $\nu = (1^{m_1}2^{m_2} \cdots )$, we establish the
notations:
\begin{eqnarray}
   \wa_{-\nu}^{\mathbb T}
&=&\prod_j (\wa_{-j}^{\mathbb T})^{m_j}, \label{wa_nu} \\
   \mathfrak z_\nu
&=&\prod_j (j^{m_j}m_j!).   \label{z_nu}
\end{eqnarray}
By Theorem~\ref{thm_struc}, a $\C$-linear basis of the space
$\Wh_{\mathbb T}$ is given by
\begin{eqnarray}   \label{WB2}
  \W {\mathcal B}_2
= \left \{ \big (\wt^{\mathbb T} \big )^{i} \,
   \wa_{-\nu}^{\mathbb T} \vac \right \}_{i \ge 0, \, \nu}.
\end{eqnarray}

\begin{lemma} \label{mon_pairing}
$\left \langle \big (\wt^{\mathbb T} \big )^{i}\wa_{-\nu}^{\mathbb
T}\vac, \big (\wt^{\mathbb T} \big )^{\w i}\wa_{-\w \nu}^{\mathbb
T}\vac \right \rangle = \mathfrak z_\nu \, \delta_{i, \w i} \,
\delta_{\nu, \w \nu}$.
\end{lemma}
\noindent {\it Proof.} By Proposition~\ref{comm_wtwa}~(i) and
Proposition~\ref{prop_hei_comm}, we obtain
\begin{eqnarray*}
&\left \langle \wt^{\mathbb T} \W A, \wt^{\mathbb T} \W B
   \right \rangle
   = \left \langle \W A, \big (\wt^{\mathbb T}
     \big )^\dagger \wt^{\mathbb T} \W B \right \rangle
   = \left \langle \W A, \W B \right \rangle,&   \\
&\left \langle \wa_{-j}^{\mathbb T} \W A, \wa_{-j}^{\mathbb T}
   \W B \right \rangle
   = \left \langle \W A, \wa_{j}^{\mathbb T}
     \wa_{-j}^{\mathbb T} \W B \right \rangle
   = j \left \langle \W A, \W B \right \rangle +
     \left \langle \W A, \wa_{-j}^{\mathbb T}
     \wa_{j}^{\mathbb T} \W B \right \rangle.&
\end{eqnarray*}
Now the lemma follows from repeatedly applying these two formulas
and
\begin{equation}
\big \langle \vac, \vac \big \rangle = \langle [C^z], [C^z] \rangle
= \langle -t, -t \rangle = 1.                     \tag*{$\qed$}
\end{equation}

\subsection{Relations with the Heisenberg operators on
             $\bigoplus_n H^{2n}_{\mathbb T}(S^{[n]})$}
\label{subsect_relation} $\,$
\par

For $n \ge 0$, let $\mathbb H_{\mathbb T, n} = H_{\mathbb
T}^{2n}(S^{[n]})$. Define the infinite dimensional space
\begin{eqnarray}      \label{fock}
\fock = \bigoplus_{n=0}^{+\infty} \mathbb H_{\mathbb T, n}.
\end{eqnarray}
In \cite{Vas} (see also \cite{LQW2}), an irreducible representation
of a Heisenberg algebra on the space $\fock$ was constructed. The
Heisenberg algebra is generated by the linear operators $\mathfrak
a_n^{\mathbb T} := \mathfrak a_n([C^z])$ in $\End(\fock)$ and the
identity operator ${\rm Id}_{\fock}$.

The operators $\mathfrak a_n([C^z])$ were defined similarly as in \S
\ref{subsect_hei}. Let $Y = C^w$ or $C^z$. For $m \ge 0$ and $n >
0$, define the closed subset $Q^{[m+n,m]}_Y$ of $\hil{m+n} \times
\hil{m}$:
\begin{eqnarray*}
Q^{[m+n,m]}_Y= \big \{ (\xi, \eta) |\, \xi\supset \eta, \,
\Supp(I_\eta/I_\xi)= \{ s \} \subset Y \big \}.
\end{eqnarray*}
Define $\mathfrak a_0([Y]) = 0 \in \End(\fock)$, and define
$\mathfrak a_{-n}([Y]), \mathfrak a_{n}([Y]) \in \End(\fock)$ by
\begin{eqnarray*}
   \mathfrak a_{-n}([Y])(A)
&=&D^{-1} (p_1')_! (\text{Id}_{\hil{m+n}} \times
     \iota_{m})_!^{-1} \left ( p_2^*A \cap
     \big [Q^{[m+n,m]}_Y \big ]\right ),              \\
   \mathfrak a_{n}([Y])(B)
&=& (-1)^n \,\, D^{-1}(p_2')_! (\iota_{m+n} \times
     \text{Id}_{\hil{m}})_!^{-1} \left (p_1^*B \cap
     \big [Q^{[m+n,m]}_Y \big ]\right )
\end{eqnarray*}
for $A \in \mathbb H_{\mathbb T, m} = H_{\mathbb T}^{2m}(S^{[m]})$
and $B \in \mathbb H_{\mathbb T, m+n}$, where $p_1', p_2'$ are the
projections:
\begin{eqnarray*}
&p_1':&\hil{m+n} \times \big ( \hil{m}
          \big )^{\mathbb T} \to \hil{m+n}, \\
&p_2':&\big ( \hil{m+n} \big )^{\mathbb T} \times
          \hil{m} \to \hil{m},
\end{eqnarray*}
$p_1, p_2$ are the two projections on $\hil{m+n} \times \hil{m}$,
and $\iota_m$ is the inclusion map:
\begin{eqnarray*}
\iota_m: \,\, (\hil{m})^{\mathbb T} \to \hil{m}.
\end{eqnarray*}
It was proved in \cite{Vas} (see \cite{LQW2} for the correct sign)
that
\begin{eqnarray}        \label{vas_hei}
[\mathfrak a_n^{\mathbb T}, \mathfrak a_m^{\mathbb T}] = n
\delta_{n, -m} \,\, \text{\rm Id}_{\fock}.
\end{eqnarray}

Recall the morphism $\rho_m: \hil{m, m+1} \to S$ from (\ref{rho_n}).
In addition, there are two natural morphisms from $\hil{m, m+1}$ to
$\hil{m}$ and $\hil{m+1}$ respectively:
\begin{eqnarray*}
\begin{array}{ccc}
\hil{m, m+1}&\overset {g_{m+1}} \longrightarrow&\hil{m+1}\\
{\,\,\,} \downarrow {f_m}&&\\
\hil{m}.&&
\end{array}
\end{eqnarray*}
These morphisms $\rho_m$, $f_m$ and $g_{m+1}$ are $\mathbb
T$-equivariant, and there are induced maps:
\begin{eqnarray}
&t \cup f_{m}^*:&\mathbb H_{\mathbb T, m} \to
                 \Wh_{\mathbb T, m},   \label{fm*}   \\
&g_{m+1}^*:&\mathbb H_{\mathbb T, m+1} \to
                 \Wh_{\mathbb T, m}.   \label{gm+1*}
\end{eqnarray}
We remark that $H^2_{\mathbb T}(S) = H^2_{\mathbb T}(\C^2) = \C t$.
Therefore, $t \cup f_{m}^*$ is essentially the only way to come up
with a map $\mathbb H_{\mathbb T, m} \to \Wh_{\mathbb T, m}$ based
on the pullback $f_{m}^*$.

\begin{proposition}  \label{wa_comm_g}
{\rm (i)} For every $n \in \Z$, there is a commutative diagram:
\begin{eqnarray}     \label{wa_comm_g.1}
\CD \mathbb H_{\mathbb T, m} @>{\mathfrak a_{-n}^{\mathbb T}}>>
    \mathbb H_{\mathbb T, m+n}  \\
@VV{t \cup f_m^*}V @VV{t \cup f_{m+n}^*}V \\
\Wh_{\mathbb T, m} @>{\wa_{-n}^{\mathbb T}}>>\Wh_{\mathbb T, m+n}.
\endCD
\end{eqnarray}

{\rm (ii)} Let $n > 0$. Then we have a commutative diagram:
\begin{eqnarray} \label{wa_comm_g.2}
\CD \mathbb H_{\mathbb T, m+1} @<{\mathfrak a_{n}^{\mathbb T}}<<
    \mathbb H_{\mathbb T, m+n+1}  \\
@VV{g_{m+1}^*}V @VV{g_{m+n+1}^*}V \\
\Wh_{\mathbb T, m} @<{\wa_{n}^{\mathbb T}}<<\Wh_{\mathbb T, m+n};
\endCD
\end{eqnarray}

{\rm (iii)} Let $n > 0$ and $A \in \mathbb H_{\mathbb T, m+1}$.
Then, we have
\begin{eqnarray}
   g_n^*\mathfrak a_{-n}^{\mathbb T}\vac
&=&n \cdot (\wt^{\mathbb T})^{n-1}\vac,     \label{wa_comm_g.3'}  \\
   g_{m+n+1}^* \mathfrak a_{-n}^{\mathbb T}(A)
&=&\wa_{-n}^{\mathbb T}(g_{m+1}^*A) - n \cdot (\wt^{\mathbb
T})^{n-1}
   \big ( t \cup f_{m+1}^*(A) \big ).   \label{wa_comm_g.3}
\end{eqnarray}
\end{proposition}
\begin{proof}
Note that $\mathfrak a_{-n}^{\mathbb T} = \mathfrak a_{-n}([C^z])$
and $\wa_{-n}^{\mathbb T} = \wa_{-n}([C^z])$. Regard the two
operators $\mathfrak a_{-n}^{\mathbb T}$ and $\wa_{-n}^{\mathbb T}$
as operators on the vector spaces
\begin{eqnarray}   \label{Wft'1}
\fock' = \bigoplus_{m=0}^{+\infty} H^*_{\mathbb T}(S^{[m]})',
              \quad
\Wft' =\bigoplus_{m=0}^{+\infty} H^*_{\mathbb T}(S^{[m, m+1]})'
\end{eqnarray}
respectively. Then similar arguments as in the proofs of the
Lemma~3.4, Lemma~4.4 and Proposition~4.6 in \cite{LQ} prove the
results. One needs to note that $\wa_{-n}^{\mathbb T}$ is
$\C[t]$-linear and the sign discrepancy in (\ref{wa_comm_g.3}) 
comes from the fact $ [C^z] = -t = \vac$.
\end{proof}
\subsection{Further properties of $t \cup f_n^*$ and $g_{n+1}^*$}
\label{subsect_further} $\,$
\par

We recall some results from \cite{Vas, LQW2} first. As in
Corollary~\ref{cor_e_c}~(i),
\begin{eqnarray*}
H_{\mathbb T}^{4n}(S^{[n]}) = t^n \cup H_{\mathbb T}^{2n}(S^{[n]}).
\end{eqnarray*}
Also, $H_{\mathbb T}^*(S^{[n]})$ is a free $\C[t]$-module. A ring
product $\star$ on $\mathbb H_{\mathbb T, n}$ is defined by
\begin{eqnarray}      \label{*}
t^{n} \cup(A \star B)= A \cup B.
\end{eqnarray}
 For $\la \vdash n$, define the class $[\la] \in
\mathbb H_{\mathbb T, n} = H_{\mathbb T}^{2n}(S^{[n]})$ by
\begin{eqnarray}  \label{normalize.1}
t^n \cup [\la] = [\xi_\la]
\end{eqnarray}
(note that our $[\la]$ differs that in \cite{Vas, LQW2} by a
scalar). Then the classes
\begin{eqnarray*}
[\la], \quad \la \vdash n
\end{eqnarray*}
form a linear basis of $\mathbb H_{\mathbb T, n}$. Let $\iota_{\la}:
\,\, \xi_{\la} \hookrightarrow S^{[n]}$ be the inclusion map, and
let
\begin{eqnarray*}
\iota_n = \bigoplus\limits_{\la \vdash n} \iota_{\la}: \,\,
(\hil{n})^{\mathbb T} \to \hil{n}.
\end{eqnarray*}
The inverse of the induced Gysin map $\iota_{n!} \colon H^*_{\mathbb
T}\big ( (S^{[n]})^{\mathbb T}\big )' \longrightarrow H^*_{\mathbb
T}(S^{[n]} )'$ is given by
\begin{eqnarray}   \label{loc_for.1}
(\iota_{n!})^{-1}(A) \,\, = \,\, \left (
\frac{(\iota_{\la})^*(A)}{(-1)^{n} h(\la)^2 \, t^{2n}} \right )_{\la
\vdash n}.
\end{eqnarray}
Define a bilinear pairing on the localization $H^*_{\mathbb
T}(S^{[n]})'$ by:
\begin{eqnarray}
&\langle -, - \rangle:
   H^*_{\mathbb T}(S^{[n]})' \otimes_{\mathbb C[t]'}
   H^*_{\mathbb T}(S^{[n]})'\to \mathbb C[t]', &
   \label{pairing.3} \\
&\langle A, B \rangle = (-1)^{n} \,
   \pi_{n!}(\iota_{n!})^{-1}(A \cup B)& \label{pairing.4}
\end{eqnarray}
where $\pi_n$ is the projection of the set $(S^{[n]})^{\mathbb T}$
of ${\mathbb T}$-fixed points to a point. This induces a $\C$-valued
bilinear pairing, again denoted by $\langle -, - \rangle$, on
$\fock$ with
\begin{eqnarray}   \label{lawla}
\langle [\la], [\w \la] \rangle = \delta_{\la, \w \la} \, h(\la)^2.
\end{eqnarray}
Then the operator $\mathfrak a_n^{\mathbb T}$ with $n > 0$ is the
adjoint of $\mathfrak a_{-n}^{\mathbb T}$, and
\begin{eqnarray}   \label{mon_pairing.2}
\left \langle \mathfrak a_{-\la}^{\mathbb T}\vac, \mathfrak a_{-\w
\la}^{\mathbb T}\vac \right \rangle = \mathfrak z_\la  \,
\delta_{\la, \w \la}
\end{eqnarray}   
where for a partition $\la = (1^{m_1}2^{m_2} \cdots)$, 
$\mathfrak a_{-\la}^{\mathbb T}$ is defined by
\begin{eqnarray}   \label{a_la}
  \mathfrak a_{-\la}^{\mathbb T}
= \prod_j (\mathfrak a_{-j}^{\mathbb T})^{m_j}.  
\end{eqnarray}

\begin{proposition}  \label{prop_fg^*}
 With the product $\star$ on $\mathbb H_{\mathbb T, n}$ and
the product $\w \star$ on $\Wh_{\mathbb T, n}$, the linear maps 
$t \cup f_n^*$ and $g_{n+1}^*$ are ring homomorphisms.
\end{proposition}
\begin{proof} 
Let $A, B \in \mathbb H_{\mathbb T, n}$. By (\ref{w*}) and
(\ref{*}), we have
\begin{eqnarray*}
t^{n+1} \cup \big ((t \cup f_n^*A) \,\, \w \star \,\,
(t \cup f_n^*B)\big ) &=& (t \cup f_n^*A) \cup (t \cup f_n^*B),\\
t^{n} \cup (A \star B) &=& A \cup B.
\end{eqnarray*}
Since $f_n^*: H^*_{\mathbb T}(\hil{n}) \to H^*_{\mathbb T}(\hil{n,
n+1})$ is a ring homomorphism with respect to the cup products and
also a $\C[t]$-module homomorphism, we obtain
\begin{eqnarray*}
   t^{n+1} \cup \big ((t \cup f_n^*A) \,\, \w \star \,\,
    (t \cup f_n^*B)\big )
&=&t^2 \cup f_n^*A \cup f_n^*B           \\
&=&t^2 \cup f_n^*(A \cup B)              \\
&=&t^2 \cup f_n^*\big (t^{n} \cup (A \star B)\big )    \\
&=&t^{n+2} \cup f_n^*(A \star B).
\end{eqnarray*}
So $(t \cup f_n^*A) \,\, \w \star \,\, (t \cup f_n^*B) = t \cup
f_n^*(A \star B)$, and $t \cup f_n^*$ is a ring homomorphism.

By a similar argument, we see that $g_{n+1}^*$ is 
a ring homomorphism.
\end{proof}

\begin{remark}
(i) We can also show that the linear map $t \cup f_n^*$ preserves 
bilinear forms, and that $\langle g_{n+1}^*A, g_{n+1}^*B \rangle 
= (n+1) \, \langle A, B \rangle$ for $A, B \in 
\mathbb H_{\mathbb T, n+1}$. Moreover,
\begin{eqnarray*}
   t \cup f_{|\la|}^*[\la]
&=&-\sum_{\mu \atop (\la, \mu) \,\, \text{\rm incidence}}
     \frac{h(\la)^2}{h(\la,  \mu)} \,\, [\la,  \mu],
     \label{lma_fg^*.01}     \\
   g_{|\mu|}^*[\mu]
&=&\sum_{\la \atop (\la, \mu) \,\, \text{\rm incidence}}
\frac{h(\mu)^2}{h(\la, \mu)}
     \,\, [\la,\mu].        \label{lma_fg^*.02}
\end{eqnarray*}

(ii) It follows from (i) that the number $h(\la, \mu)$ satisfies 
some interesting identities:
\begin{eqnarray*}
   \sum_{\mu \atop (\la, \mu) \,\, \text{\rm incidence}}
      \frac{h(\la)^2}{h(\la, \mu)}
&=&\sum_{k = 0}^{s(\la)} \prod_{0 \le i \le s(\la) \atop i \ne k}
      \frac{h(\square_{k, i}')}{1+h(\square_{k, i})}
      = 1,     \\
   \sum_{\la \atop (\la, \mu) \,\, \text{\rm incidence}}
   \frac{h(\mu)^2}{h(\la, \mu)}
&=&|\mu|.
\end{eqnarray*}
\end{remark}

\section{\bf Transformations among various linear bases of $\Wh_{\mathbb T}$}
\label{sect_Transf}

Recall that the vector space $\Wh_{\mathbb T}$ has two linear bases:
\begin{eqnarray*}
   \W {\mathcal B}_1
    = \big \{ [\la, \mu] \big \}_{(\la, \mu) \,\,
          \text{incidence}},\quad
   \W {\mathcal B}_2
=\left \{ \big (\wt^{\mathbb T} \big )^{i} \,
   \wa_{-\nu}^{\mathbb T} \vac \right \}_{i \ge 0, \, \nu}.
\end{eqnarray*}
The ring structure of $\Wh_{\mathbb T, n}$ is easily described in
terms of the linear basis $\W {\mathcal B}_1$ coming from the fixed points.
However, the basis $\W {\mathcal B}_1$ doesn't exist in the ordinary
cohomology of $\hil{n,n+1}$, while the second basis $\W {\mathcal
B}_2$ survives. Since we are interested in the ring structure of 
the ordinary cohomology $H^*(\hil{n,n+1})$, it is important to know 
the linear transformation between these two bases. In this section, 
we give an algorithm to express $\big (\wt^{\mathbb T} \big )^{i} \, 
\wa_{-\nu}^{\mathbb T} \vac$ as a linear combination of the elements
in $\W {\mathcal B}_1$. This allows us to describe (implicitly) the
ring structure of $\Wh_{\mathbb T, n}$ in terms of the  linear basis
$\W {\mathcal B}_2$. The method is to introduce the subvarieties $\W
L^{\la, \mu}(C^z)$ on $\hil{n,n+1}$ and study the actions of the loop
algebra and the Heisenberg algebra on them. As a consequence, 
the classes $[\W L^{\la, \mu}(C^z)]$ provide a link between 
the bases $\W {\mathcal B}_1$ and $\W {\mathcal B}_2$.

\subsection{The subvariety $L^\la(C^z)$ on $\hil{n}$}
\label{subsect_LlaC} $\,$
\par

Let $C = C^z$ be the $z$-axis of $S = \C^2$. For $\la = (\la_1 \ge
\la_2 \ge \ldots ) \vdash n$, let
\begin{eqnarray*}
S^n_\la C = \left \{ \sum_i \la_i s_i \in {\rm Sym}^n(S) | s_i \in C
\,\, \text{\rm and the $s_i$'s are distinct} \right \}.
\end{eqnarray*}
Recall the Hilbert-Chow morphism $\pi_n: \Sn \rightarrow {\rm
Sym}^n(S)$. Let
\begin{eqnarray}  \label{LlaC}
L^\la C = \text{\rm Closure of } (\pi_n)^{-1}(S^n_\la C).
\end{eqnarray}
The subvariety $L^\la C$ was first introduced in \cite{Gro}, and was
studied intensively in \cite{Na2, Na3}. Note that $L^\la C$ is
irreducible, of dimension $n$ and $\mathbb T$-invariant.

By the results in \cite{Na2, Na3, Vas}, $\xi_\la$ is a smooth point
of $L^\la C$, $\xi_{\w \la} \in L^\la C$ if and only if $\w \la \le
\la$ where $\le$ denotes the dominance partial ordering, and
\begin{eqnarray}   \label{aiLlaC.0}
\mathfrak a_{-m}^{\mathbb T} [L^\lambda C] = \sum_{\nu} a_{\lambda,
\nu} [L^{\nu} C],
\end{eqnarray}
where the summation is over partitions $\nu$ of $|\lambda|+m$, which
are obtained as follows:
\begin{enumerate}
\item[{\rm (i)}] add $m$ to a term in $\lambda$,
say $\lambda_k$ (possibly $0$), and then
\item[{\rm (ii)}] arrange it in descending order.
\end{enumerate}
The coefficient $a_{\lambda, \nu}$ is the number of $\ell$ with
$\nu_\ell=\lambda_k+m$. Denote the above $\nu$ by:
\begin{eqnarray}   \label{wla}
\nu = \la(\la_k, m).
\end{eqnarray}
Then, formula (\ref{aiLlaC.0}) can be rewritten as
\begin{eqnarray}   \label{aiLlaC}
\mathfrak a_{-m}^{\mathbb T} [L^\lambda C] = \sum_{{\rm
distinct}\,\, \la_k} a_{\lambda, \la(\la_k, m)} 
[L^{\la(\la_k, m)} C].
\end{eqnarray}
\subsection{The subvariety $\W L^{\la, \mu}(C^z)$ on $\hil{n,n+1}$}
\label{subsect_LlamuC} $\,$
\par

 Again, let $C = C^z$ be the $z$-axis of $S =
\C^2$. Let $(\la, \mu)$ be an incidence pair:
\begin{eqnarray*}
\la &=& (\cdots (i-1)^{m_{i-1}} i^{m_i}
          (i+1)^{m_{i+1}} \cdots) \vdash n, \\
\mu &=& (\cdots (i-1)^{m_{i-1}} i^{m_i-1}
          (i+1)^{m_{i+1}+1} (i+2)^{m_{i+2}} \cdots)
\end{eqnarray*}
where the parts $(i-1)^{m_{i-1}}$, $i^{m_i}$ and $i^{m_i-1}$ do not
appear when $i = 0$. We fix this $i$ throughout the subsection. Let
$\W L^{\la, \mu}C$ be the closed subset of $\hil{n, n+1}$ defined by
\begin{eqnarray}     \label{LlamuC}
\W L^{\la, \mu}C = \left \{ (\xi, \xi^\prime)|\xi\in L^\la C,\,\,
\xi^\prime\in L^\mu C,\, \, \xi \subset \xi^\prime \right \}.
\end{eqnarray}
Then, $\W L^{\la, \mu}C$ is irreducible, of dimension $(n+1)$ and
$\mathbb T$-invariant.

\begin{lemma}   \label{wtLlamu}
$\wt^{\mathbb T} [\W L^{\la,\mu} C] = [\W L^{\mu, \nu} C]$ where the
partition $\nu$ is defined by
\begin{eqnarray*}
\nu = (\cdots (i-1)^{m_{i-1}}i^{m_i-1}(i+1)^{m_{i+1}}
        (i+2)^{m_{i+2}+1}(i+3)^{m_{i+3}}\cdots).
\end{eqnarray*}
\end{lemma}
\noindent
{\it Proof.} 
Let $\tilde{p}_1, \tilde{p}_2$ be the two projections
of $S^{[n+1, n+2]} \times S^{[n,n+1]}$. Then,
\begin{eqnarray*}
\wt^{\mathbb T} [\W L^{\la,\mu} C] = D^{-1}\tilde{p}_{1*} \left (
\tilde{p}_2^*[\W L^{\la,\mu} C] \cap [\Wq_n] \right ).
\end{eqnarray*}
An element $\big ( (\xi^{\prime}, \xi^{\prime \prime}), (\xi,
\xi^{\prime}) \big )$ in $\Wq_n \cap (\w p_2)^{-1} \W L^{\la,\mu}C$
is of the form
\begin{eqnarray*}
\xi=\xi_0+\xi_s \in L^\la C, \quad \xi^{\prime}=\xi_0+\xi^\prime_s
\in L^\mu C, \quad \xi^{\prime\prime}=\xi_0+\xi^{\prime\prime}_s
\end{eqnarray*}
where $\xi_s\subset\xi^\prime_s\subset \xi^{\prime\prime}_s$,
$\Supp(\xi^{\prime\prime}_s) = \{s\} \subset C$, and $s \notin
\Supp(\xi_0) \subset C$. The image $p_1 \left (\Wq_n \cap (\w
p_2)^{-1} \W L^{\la,\mu}C \right )$ consists of elements of the form
\begin{eqnarray*}
\xi^{\prime}=\xi_0+\xi^\prime_s,\quad
\xi^{\prime\prime}=\xi_0+\xi^{\prime\prime}_s.
\end{eqnarray*}
Let $\ell(\xi_s) = \ell$. Choose a local coordinate $(w_s, z)$ of
$S$ near $s$ such that $C$ is given locally near $s$ by $w_s=0$. Now
there are two cases:

\medskip\noindent
{\bf Case 1:} the element $\xi^{\prime\prime}_s \in M_{\ell+2}(s)$
is not generic. Since $M_{\ell+1, \ell+2}(s)$ is irreducible and has
dimension $(\ell + 1)$, the corresponding element $(\xi^{\prime},
\xi^{\prime\prime}) = (\xi_0+\xi^\prime_s,
\xi_0+\xi^{\prime\prime}_s)$ forms a subset of dimension less than
$(n+2)$.

\medskip\noindent
{\bf Case 2:} the element $\xi^{\prime\prime}_s \in M_{\ell+2}(s)$
is generic. Then $\xi^{\prime\prime}_s$ is of the form:
\begin{eqnarray*}
I_{\xi^{\prime\prime}_s} = (w_s^{\ell+2}, z+b_1w_s+\ldots+b_{\ell+1}
w_s^{\ell+1})
\end{eqnarray*}
where $b_1, \ldots, b_{\ell+1} \in \C$. The corresponding
$\xi^\prime_s$ and $\xi_s$ must be of the form
\begin{eqnarray*}
I_{\xi^{\prime}_s} &=& (w_s^{\ell+1},
  z+b_1w_s+\ldots+b_{\ell} w_s^{\ell}),  \\
I_{\xi_s} &=& (w_s^{\ell},
  z+b_1w_s+\ldots+b_{\ell-1} w_s^{\ell-1}).
\end{eqnarray*}
By the definition of $L^\mu C$, we must have $\ell = i$.

It follows that $p_1 \left (\Wq_n \cap (\w p_2)^{-1} \W L^{\la,\mu}C
\right )$ has exactly one irreducible component of the expected
dimension $(n+2)$, which is $\W L^{\mu, \nu}C$, and possibly other
components with smaller dimension. Note that the intersection $\Wq_n
\cap (\w p_2)^{-1} \W L^{\la,\mu}C$ is transversal along the element
$\big ( (\xi^{\prime}, \xi^{\prime \prime}), (\xi, \xi^{\prime})
\big )$ if $\xi^{\prime\prime}_s$ is from Case 2 and $\xi_0$ is
generic. Hence
\begin{equation}
\wt^{\mathbb T} [\W L^{\la,\mu} C] = [\W L^{\mu, \nu} C].
\tag*{$\qed$}
\end{equation}

\begin{lemma}   \label{amLlamu}
Let $m_0 = 2$. Let $\la(\la_k, m)$ and $a_{\lambda, \la(\la_k, m)}$
be from (\ref{aiLlaC}). Then,
\begin{eqnarray*}
   \wa_{-m}^{\mathbb T} [\W L^{\la, \mu} C]
&=&\sum_{{\rm distinct}\,\, \la_k \ne i} a_{\lambda, \la(\la_k, m)}
    [\W L^{\la(\la_k, m), \mu(\la_k, m)} C]  \\
& &+ \,\, (1- \delta_{m_i-1, 0}) a_{\lambda, \la(i, m)}
    [\W L^{\la(i, m), \mu(i, m)} C] + [\W L^{\la(i, m), \mu(i+1, m)} C].
\end{eqnarray*}
\end{lemma}
\begin{proof}
Let $n = |\la|$. Recall that $\wa_{-m}^{\mathbb T} = \wa_{-m}([C])$.
By (\ref{waY}), we have
\begin{eqnarray*}
\wa_{-m}^{\mathbb T}([C]) [\W L^{\la,\mu} C] = D^{-1}\tilde{p}_{1*}
\left ( \tilde{p}_2^*[\W L^{\la,\mu} C] \cap [\Wq_C^{[n+m, n]}]
\right ).
\end{eqnarray*}

Let $\xi_{0, 1} = \xi_{0, 2} = \emptyset$. Then a generic element
$(\xi, \xi^\prime) \in \W L^{\la,\mu} C$ is of the form:
\begin{eqnarray*}
\xi &=& \sum_{r \ge 0} \sum_{1 \le j \le m_r} \xi_{r, j}, \\
\xi^\prime &=& \xi - \xi_{i, m_i} + \xi_{i, m_i}^\prime
\end{eqnarray*}
where $i \ge 0$, $\ell(\xi_{r, j}) = r$, $\Supp(\xi_{r, j}) =
\{s_{r, j}\} \subset C$ for $r \ge 1$, the points $s_{r, j}$ are
distinct, $(\xi_{i, m_i}, \xi_{i, m_i}^\prime) \in M_{i, i+1}(s_{i,
m_i})$ when $i > 0$, and when $i = 0$, $\xi_{i, m_i}^\prime$ is a
point in $C$ different from the points $s_{r, j}$ with $r \ge 1$.

The effect of the action of $\wa_{-m}^{\mathbb T}([C])$ on $(\xi,
\xi^\prime) \in \W L^{\la,\mu} C$ has two types:

\medskip\noindent
{\bf Type 1}: the action results in generic elements $(\eta,
\eta^\prime)$ of the form:
\begin{eqnarray*}
\eta &=& \xi - \xi_{r_0, j_0} + \eta_{r_0 + m}, \\
\eta^\prime &=& \xi^\prime - \xi_{r_0, j_0} + \eta_{r_0 + m}
\end{eqnarray*}
where $(r_0, j_0) \ne (i, m_i)$, $\ell(\eta_{r_0 + m}) = r_0 + m$,
$\Supp(\eta_{r_0 + m}) = \{ s_{r_0 + m}\} \subset C$, and
\begin{eqnarray*}
s_{r_0 + m} \not \in \Supp(\xi^\prime).
\end{eqnarray*}
This type of action is similar to the action of $\mathfrak
a_{-m}^{\mathbb T}$ on $[L^\lambda C]$. It follows that
$\wa_{-m}^{\mathbb T} [\W L^{\la, \mu} C]$ contains $a_{\lambda,
\la(\la_k, m)} [\W L^{\la(\la_k, m), \mu(\la_k, m)} C]$ when $\la_k
\ne i$, and contains
\begin{eqnarray*}
(1- \delta_{m_i-1, 0}) a_{\lambda, \la(i, m)} [\W L^{\la(i, m),
\mu(i, m)} C]
\end{eqnarray*}
which should be regarded as nonexistent when $m_i = 1$.

\medskip\noindent
{\bf Type 2}: the action results in generic elements $(\eta,
\eta^\prime)$ of the form:
\begin{eqnarray*}
\eta &=& \xi - \xi_{i, m_i} + \eta_{i + m}, \\
\eta^\prime &=& \xi - \xi_{i, m_i} + \eta_{i + m}^\prime
\end{eqnarray*}
where $i \ge 0$, $(\eta_{i + m}, \eta_{i + m}^\prime) \in M_{i + m,
i+m+1}(s_{i + m})$ for some $s_{i + m} \in C$ with
\begin{eqnarray*}
s_{i + m} \not \in \Supp(\xi).
\end{eqnarray*}
It follows that $\wa_{-m}^{\mathbb T} [\W L^{\la, \mu} C]$ contains
$a \, [\W L^{\la(i, m), \mu(i+1, m)} C]$ for certain multiplicity
$a$. As at the end of the proof of Lemma~\ref{wtLlamu}, the
intersection multiplicity $a$ is $1$.
\end{proof}

\begin{proposition}   \label{Llamu_mono}
Let $(\la, \mu)$ be an incidence pair. Then, there exists an
algorithm to express $[\W L^{\la, \mu} C]$ as a linear combination
of the elements in the linear basis $\W {\mathcal B}_2$.
\end{proposition}
\begin{proof}
Use induction on $|\la|$. When $|\la| = 0$, we have
\begin{eqnarray}        \label{Llamu_mono.1}
[\W L^{\la, \mu} C] = [C] = [C^z] = \vac.
\end{eqnarray}
So the conclusion holds. Next, let $n \ge 1$ and assume that the
conclusion holds for $[\W L^{\la, \mu} C]$ whenever $|\la| < n$. In
the following, let $|\la| = n$.

Use a second induction on $\ell(\la)$. When $\ell(\la) = 1$, $\la =
(n)$, and either $\mu = (1, n)$ or $\mu = (n+1)$. When $\mu =
(n+1)$, we see from Lemma~\ref{wtLlamu} that
\begin{eqnarray*}
[\W L^{\la, \mu} C] = [\W L^{(n), (n+1)} C] = (\wt^{\mathbb T})^n
\vac,
\end{eqnarray*}
and so the conclusion holds in this case. When $\mu = (1, n)$,
applying the operator $\wa_{-n}^{\mathbb T}$ to (\ref{Llamu_mono.1})
and using Lemma~\ref{amLlamu}, we conclude that
\begin{eqnarray*}
[\W L^{(n), (n+1)} C] + [\W L^{(n), (1, n)} C] = \wa_{-n}^{\mathbb
T} \vac.
\end{eqnarray*}
Hence, the conclusion holds for $[\W L^{(n), (1, n)} C]$ as well.

Let $|\la| = n$ and $\ell(\la) > 1$. In this case, we can choose a
part $m > 0$ of $\la$ so that $m$ is also a part of $\mu$. Let
$\la^{[m]}$ and $\mu^{[m]}$ be the partitions of $n - m$ obtained
from $\la$ and $\mu$ respectively by deleting a copy of part $m$.
Let
\begin{eqnarray*}
\la &=& \big (\cdots (i-1)^{m_{i-1}} i^{m_i}
          (i+1)^{m_{i+1}} \cdots \big ) \vdash n, \\
\mu &=& \big (\cdots (i-1)^{m_{i-1}} i^{m_i-1}
          (i+1)^{m_{i+1}+1} (i+2)^{m_{i+2}} \cdots \big ).
\end{eqnarray*}
Apply the formula in Lemma~\ref{amLlamu} to 
$\wa_{-m}^{\mathbb T} [\W L^{\la^{[m]}, \mu^{[m]}} C]$, 
for which the conclusion holds by the induction
hypothesis on the size of partitions since $|\la^{[m]}|<|\la|$. 
There are three types of terms on the right hand side of 
the formula. One is $[\W L^{\la, \mu} C]$ with a non-zero 
coefficient. The second type is $[\W L^{\w \la, \w \mu} C]$ 
with $|\w \la| = |\la| = n$ and $\ell(\w\la )<\ell(\la)$, 
for which the conclusion
holds by the induction hypothesis on the length of partitions. 
The third type is the term $[\W L^{\la^{[m]}(i, m), 
\mu^{[m]}(i+1, m)} C]=(\wt^{\mathbb T})^m[\W
L^{\la^{[m]}, \mu^{[m]}} C]$, for which the conclusion holds by
induction on the size again since $|\la^{[m]}|<|\la|$. 
It follows immediately that the conclusion holds for 
$[\W L^{\la, \mu} C]$.
\end{proof}
\subsection{Transformations between the linear bases
         $\W {\mathcal B}_1$ and $\W {\mathcal B}_2$}
\label{subsect_B1B2} $\,$
\par

Note that the $\mathbb T$-action on $S^{[n, n+1]}$ induces a cell
decomposition of
\begin{eqnarray*}
\W L^{n, n+1} C := \coprod_{\la \vdash n \atop (\la, \mu) \,\, {\rm
incidence}} \W L^{\la, \mu} C.
\end{eqnarray*}
Let $(\xi, \xi') \in \W L^{\la, \mu} C$ be a generic point. Then we
see that
\begin{eqnarray}   \label{T_cell}
\lim_{a \to 0} a (\xi, \xi') = (\xi_\la, \xi_\mu)
\end{eqnarray}
for $a \in \mathbb T$. Since $\dim \W L^{\la, \mu} C = \dim \W L^{n,
n+1} C = n+1$, we conclude that the cell $\mathcal C^{\la, \mu}$
corresponding to the fixed point $(\xi_\la, \xi_\mu)$ is isomorphic
to $\C^{n+1}$ and $\W L^{\la, \mu} C$ is the closure of $\mathcal
C^{\la, \mu}$. In particular, the fixed point $(\xi_\la, \xi_\mu)$
is a smooth point of $\W L^{\la, \mu} C$.

Recall the notation $k(\la, \mu)$ from
Definition~\ref{squareii'}~(ii). Put
\begin{eqnarray}   \label{h'lagamma}
h_+(\la, \mu) = h(\la) \cdot \prod_{k(\la, \mu) + 1 \le i \le
s(\la)} \frac{1+h(\square_{k(\la, \mu), i})} {h(\square_{k(\la,
\mu), i}')}.
\end{eqnarray}

\begin{lemma}   \label{tang_Llamu}
Let $(\xi_{\la}, \xi_{\mu}) \in \hil{n, n+1}$ be a $\mathbb T$-fixed
point. Then the $\mathbb T$-equivariant Euler class 
of the tangent space of $\W L^{\la, \mu} C$ at $(\xi_{\la},
\xi_{\mu})$ is equal to $h_+(\la, \mu) t^{n+1}$.
\end{lemma}
\begin{proof}
Since $\W L^{\la, \mu} C$ is the closure of the cell $\mathcal
C^{\la, \mu}$ corresponding to the fixed point $(\xi_\la, \xi_\mu)$,
it suffice to compute the $\mathbb T$-equivariant Euler class of the
tangent space of $\mathcal C^{\la, \mu}$ at $(\xi_{\la},
\xi_{\mu})$. By (\ref{T_cell}), the tangent space of 
$\mathcal C^{\la, \mu}$ at $(\xi_{\la}, \xi_{\mu})$ 
is the positive part of the tangent space
of $\hil{n, n+1}$ at $(\xi_{\la}, \xi_{\mu})$. The positive part of
the tangent space of $\hil{n, n+1}$ at $(\xi_{\la}, \xi_{\mu})$ can
be read from the detailed study of the equivariant Zariski tangent
space in the Appendix. Hence we see that the $\mathbb
T$-equivariant Euler class of the tangent space of $\mathcal C^{\la,
\mu}$ (and hence of $\W L^{\la, \mu} C$) at $(\xi_{\la}, \xi_{\mu})$
is $h_+(\la, \mu) t^{n+1}$.
\end{proof}

Note that $(\xi_{\w \la}, \xi_{\w \mu}) \in \W L^{\la, \mu} C$ only
if $\la \ge \w \la$ and $\mu \ge \w \mu$. When $\la \ge \w \la$,
$\mu \ge \w \mu$ and $(\la, \mu) \ne (\w \la, \w \mu)$, we define
$(\la, \mu) > (\w \la, \w \mu)$.

\begin{lemma}   \label{Llamu_fixed}
Let $(\la, \mu)$ be an incidence pair. Then, we have
\begin{eqnarray}   \label{Llamu_fixed.0}
[\W L^{\la, \mu} C] = h_+(\la, \mu)^{-1} \, [\la, \mu] + \sum_{(\la,
\mu) > (\w \la, \w \mu)} d_{(\la, \mu), (\w \la, \w \mu)} [\w \la,
\w \mu]
\end{eqnarray}
for some constants $d_{(\la, \mu), (\w \la, \w \mu)} \in \mathbb Q$.
Moreover, there exists an algorithm to compute all the constants
$d_{(\la, \mu), (\w \la, \w \mu)} \in \mathbb Q$.
\end{lemma}
\begin{proof}
It is known from \cite{Bri} that if $X, Y$ are $\mathbb
T$-equivariant equidimensional varieties such that $Y \subset X$ is
closed and $X^{\mathbb T}$ is finite, then
\begin{eqnarray*}
[Y] = \sum_{y \in Y^{\mathbb T}} c_y(Y) t^{- \dim Y} \cup [y] \in
H^*_{\mathbb T}(X)'
\end{eqnarray*}
for some constants $c_y(Y) \in \mathbb Q$. Moreover, if $y \in
Y^{\mathbb T}$ is a smooth point of $Y$, then $c_y(Y) \ne 0$ and
$c_y(Y)^{-1} t^{\dim Y}$ is the $\mathbb T$-equivariant Euler class
of the tangent space of $Y$ at $y$. Apply this formula to $X =
\hil{|\la|, |\la|+1}$ and $Y = \W L^{\la, \mu} C$. We see that
(\ref{Llamu_fixed.0}) follows from Lemma~\ref{tang_Llamu} and the
definition of the class $[\la, \mu]$.

By Proposition~\ref{Llamu_mono} and Lemma~\ref{mon_pairing}, there
exists an algorithm to compute the pairings among the classes $[\W
L^{\la, \mu} C]$. By (\ref{Llamu_fixed.0}), the classes $[\W L^{\la,
\mu} C]$ are related to the classes $[\la, \mu]$ via an upper
triangular matrix. The diagonal entries of the matrix are
$h_+(\la,\mu)^{-1}$, which are nonzero, and the entries above the
diagonal  are the constants $d_{(\la, \mu), (\w \la, \w \mu)}$.
Since the pairings between the fixed point classes $[\la, \mu]$  are
already computed in (\ref{muinuj}), this upper triangular
matrix can be determined. Therefore, there exists an algorithm to
compute the constants $d_{(\la, \mu), (\w \la, \w \mu)}$.
\end{proof}

\begin{theorem}  \label{thm_alg}
There exists an algorithm to express each element $\big
(\wt^{\mathbb T} \big )^{i} \, \wa_{-\nu}^{\mathbb T} \vac$ in the
linear basis $\W {\mathcal B}_2$ as a linear combination of the
elements in the linear basis $\W {\mathcal B}_1$.
\end{theorem}
\begin{proof}
By Lemma~\ref{Llamu_fixed}, there exists an algorithm to express
each element $[\la, \mu]$ in $\W {\mathcal B}_1$ as a linear
combination of the elements in the third linear basis
\begin{eqnarray}   \label{thm_alg.1}
\W {\mathcal B}_3 = \left \{ [\W L^{\la, \mu} C] \right \}_{(\la,
\mu) \,\, {\rm incidence}}
\end{eqnarray}
of $\Wh_{\mathbb T}$. Note that the transition matrix between these
two linear bases can be arranged to be lower triangular. By
Proposition~\ref{Llamu_mono}, there exists an algorithm to express
each element $[\la, \mu]$ in $\W {\mathcal B}_1$ as a linear
combination of the elements in $\W {\mathcal B}_2$. Therefore, we
conclude that there exists an algorithm to express each element
$\big (\wt^{\mathbb T} \big )^{i} \, \wa_{-\nu}^{\mathbb T} \vac$ in
$\W {\mathcal B}_2$ as a linear combination of the elements in $\W
{\mathcal B}_1$.
\end{proof}

\section{\bf Applications}
\label{sect_appl}
\subsection{Application to the ring structure of $H^*(\hil{n, n+1})$}
\label{subsect_app_ring} $\,$
\par

In \S 5.2 of \cite{LQ}, we showed that the infinite dimensional
space
\begin{eqnarray}      \label{Wfock}
\Wfock = \bigoplus_{n=0}^{+\infty} H^*(\hil{n, n+1})
\end{eqnarray}
is a representation of the Lie algebra $\C[u^{-1}]
\otimes_\C {\w {\mathfrak h}}_S$ with a highest weight vector being
the vacuum vector
\begin{eqnarray*}
\vac = 1_S \in H^0(\hil{0,1}) = H^0(S).
\end{eqnarray*}
Here, $1_S \in H^0(S)$ is the fundamental cohomology class of $S =
\C^2$, $u^{-1}$ acts via a translation operator $\wt$ similarly
defined as in \S \ref{subsect_trans}, and $\w {\mathfrak h}_S$ is
the Heisenberg algebra generated by ${\rm Id}_{\Wfock}$ and the
Heisenberg operators $\wa_n$ with $n \in \Z$. When $n > 0$, the
creation Heisenberg operator $\wa_{-n} = \wa_{-n}(1_S)$ is defined
similarly as in \S \ref{subsect_hei}. In particular, the elements in
$\Wfock$ are of the form:
\begin{eqnarray}   \label{elements}
\wt^{i} \wa_{-n_1}^{i_1} \cdots \wa_{-n_k}^{i_k} \vac
\end{eqnarray}
where $k \ge 0$, $i \ge 0$, $i_1, \ldots, i_k > 0$, 
and $n_1, \ldots, n_k > 0$.

For a partition $\nu = (1^{m_1}2^{m_2} \cdots )$, 
we introduce the notation:
\begin{eqnarray*}
\wa_{-\nu} &=& \prod_j \wa_{-j}^{m_j}.
\end{eqnarray*}

Recall the ring isomorphism in Corollary~\ref{cor_e_c}~(ii) induced
by the forgetful map
\begin{eqnarray*}
\Psi: H_{\mathbb T}^*(S^{[n, n+1]}) \to H^*(S^{[n, n+1]}).
\end{eqnarray*}
Let $\nu$ be a partition with $i + |\nu| = n$. As in \S 4 of
\cite{LQW3}, we have
\begin{eqnarray}   \label{Psi}
\Psi \left ( (-t)^{-\ell(\nu) -1} \,\, \big (\wt^{\mathbb T} \big
)^{i} \wa_{-\nu}^{\mathbb T} \vac \right ) = \wt^{i} \wa_{-\nu} \vac
\end{eqnarray}
noting $[\C^2] = -t^{-1} [C^z]$ in $H^*_{\mathbb T}(\C^2) =
H^*_{\mathbb T}(S)$, $\wa_n^{\mathbb T} = \wa_n([C^z])$, and the
Heisenberg commutation relation (\ref{thm_struc}). It follows that
the cup products of the classes $\wt^{i} \wa_{-\nu} \vac$ can be
reduced,  by using (\ref{Psi}) and Theorem~\ref{thm_alg}, to
computations in terms of the linear basis $\W {\mathcal B}_1$ of
the fixed points. The cup products of the elements in 
$\W {\mathcal B}_1$ are already determined in \S \ref{subsect_wht}. 
This gives the ring structure of $H^*(\hil{n, n+1})$.

\subsection{Application to the ring of symmetric functions}
\label{subsect_app_symme} $\,$
\par

Let $\Lambda$ be the space of symmetric functions in
infinitely many variables (see p.19 of \cite{Mac}). For a partition
$\la$, let $p_\la, m_\la$ and $s_\la$ be the power-sum symmetric
function, the monomial symmetric function and the Schur function
associated to $\la$ respectively. 
Define a ring structure on $\Lambda$ by requiring
$s_\la\cdot s_\mu=\delta_{\la, \mu}h(\la)s_\la$ for the Schur functions
$s_\la$ and $s_\mu$. Note that $\Lambda$ already has a natural ring
structure of the multiplication of functions. To avoid the
confusion, we use $(\Lambda, \cdot)$ to denote $\Lambda$ with the
new ring structure, and always refer to this new ring structure when
we mention the ring structure of $\Lambda\otimes_\Z\C$.

Let $C = C^z \subset S = \C^2$. By the Proposition~B of \cite{Vas},
there is a ring isomorphism
\begin{eqnarray}  \label{Phi}
\Phi: \mathbb H_{\mathbb T} \to (\Lambda \otimes_\Z \C, \cdot)
\end{eqnarray}
which sends the classes $\mathfrak a_{-\la}^{\mathbb T} \vac$,
$[L^\la C]$ and $(-1)^n/h(\la) \cdot [\la]$ to the symmetric
functions $p_\la, m_\la$ and $s_\la$ respectively. Moreover, under
this isomorphism, the operator $\mathfrak a_{-n}^{\mathbb T}$ 
with $n > 0$ on $\mathbb H_{\mathbb T}$ corresponds to 
multiplication by $p_{(n)}$ on $\Lambda \otimes_\Z \C$.

Extending the map $\Phi$, we define a linear isomorphism
\begin{eqnarray}  \label{WPhi}
\W \Phi: \,\, \Wh_{\mathbb T} \to \Lambda \otimes_\Z \C[v]
\end{eqnarray}
by sending $\big (\wt^{\mathbb T} \big )^{i} \, \wa_{-\la}^{\mathbb
T} \vac$ to $p_\la \otimes v^i$. Under this linear isomorphism, the
operators $\wt^{\mathbb T}$ and $\wa_{-n}^{\mathbb T}$ with $n > 0$ 
on $\Wh_{\mathbb T}$ correspond to multiplications by $1 \otimes v$ 
and $p_{(n)} \otimes 1$ on $\Lambda \otimes_\Z \C[v]$ respectively.
Moreover, the ring structure on 
$\Wh_{\mathbb T}$ induces a ring structure on 
$\Lambda \otimes_\Z \C[v]$ such that $\Lambda \otimes_\Z \C \subset 
\Lambda \otimes_\Z \C[v]$ is a subring of
$\Lambda \otimes_\Z \C[v]$. Let 
\begin{eqnarray*}
\iota \colon \Lambda\otimes_{\mathbb Z}\mathbb C
\hookrightarrow \Lambda \otimes_\Z \C[v]
\end{eqnarray*}
be the inclusion map. Next, recall from Proposition~\ref{prop_fg^*}
that there exists a ring homomorphism ${t \cup f_m^*}\colon 
\mathbb H_{\mathbb T, m}\to \Wh_{\mathbb T, m}$. It induces 
a ring homomorphism:
\begin{eqnarray*}
{t \cup f_.^*}\colon \mathbb H_{\mathbb T} \to \Wh_{\mathbb T}.
\end{eqnarray*}
By Proposition~\ref{wa_comm_g}~(i), we obtain 
a commutative diagram of ring homomorphisms:
\begin{eqnarray*}
\CD \mathbb H_{\mathbb T} @>{\Phi}>>
    \Lambda\otimes_{\mathbb Z}\mathbb C  \\
@VV{t \cup f^*_.}V @VV{\iota}V \\
\Wh_{\mathbb T} @>{\W \Phi}>>\Lambda \otimes_\Z \C[v],
\endCD
\end{eqnarray*}
which respects the Heisenberg algebra actions on 
$\mathbb H_{\mathbb T}$ and $\Wh_{\mathbb T}$.

It is natural to ask what the induced ring structure on
$\Lambda \otimes_\Z \C[v]$ is. It should provide an interesting
feature on $\Lambda \otimes_\Z \C[v]$ in the realm of symmetric
functions.

\section{\bf Appendix: the proof of Lemma~\ref{lemma_eT}}
\label{sec:appendix}

For simplicity, we denote $\square_{k(\la, \mu), i}=\square_{k, i}$
and $\square_{k(\la, \mu), i}'=\square_{k, i}'$ by $\square_{i}$ and
$\square_{i}'$ respectively. Note from (\ref{ms(la)}) and
Definition~\ref{squareii'}~(ii) that (\ref{lemma_eT.0}) is the same
as
\begin{eqnarray}   \label{lemma_eT.1}
e_{\mathbb T} =(-1)^{n+1} \cdot h(\la)^2 \cdot \prod_{0 \le i \le m
\atop i \ne k} \frac{1+h(\square_i)}{h(\square_i')} \cdot
t^{2(n+1)}.
\end{eqnarray}
To prove (\ref{lemma_eT.1}), we follow the setup in \cite{Ch1}.
There are four separate cases.

\medskip\noindent
{\bf Case 1a:} $q_k = 1$ but $p_k \ne 1$. Then $s=m$, $\alpha_k'
\not\in A$, and $\alpha_i' = \alpha_i$ if $0 \le i \le s$ and $i \ne
k$.
For $0 \le i \le (k-2)$, we have $p_i' = p_i$ and let $\beta_i \in
B$ be the $p_i$-th cell directly above $\alpha_k$. Then, $\beta_i
\in P_{\alpha_i}$. For $k \le i \le (s-1)$, we have $q_{i+1}' =
q_{i+1}$ and let $\beta_{i+1} \in B$ be the $q_{i+1}$-th cell
directly to the left of $\alpha_k$. Then, $\beta_{i+1} \in
Q_{\alpha_{i+1}}$. By the formula (2.6.1) in \cite{Ch1}, $\Hom \big
(I_{\xi_{\mu}}, R/I_{\xi_{\la}} \big )$ is equal to
\begin{eqnarray}   \label{case1a}
{\rm im}(\psi) \bigoplus \left ( \bigoplus_{i = 0}^{k-2}
     \C \phi(f_{\alpha_i, \beta_i}) \right )
\bigoplus \left ( \bigoplus_{i = k}^{s-1}
     \C \phi(f_{\alpha_{i+1}, \beta_{i+1}}) \right ).
\end{eqnarray}
Combining this with (\ref{tang.3}) and (\ref{fai'ak}), we obtain an
exact sequence
\begin{eqnarray}    \label{case1a_ker.1}
0 \to \frac{\ker(\phi-\psi)}{\bigoplus_{i=0}^s
  \C f_{\alpha_i', \alpha_k}} \to \Hom \big (I_{\xi_{\la}},
  R/I_{\xi_{\la}} \big ) \oplus {\rm im}(\psi)
  \qquad\qquad   \nonumber   \\
\to {\rm im}(\psi) \bigoplus \left ( \bigoplus_{i = 0}^{k-2}
  \C \phi(f_{\alpha_i, \beta_i}) \right )
  \bigoplus \left ( \bigoplus_{i = k}^{s-1}
  \C \phi(f_{\alpha_{i+1}, \beta_{i+1}}) \right ) \to 0.
\end{eqnarray}

If $0 \le i \le (k-1)$, then the weight of $f_{\alpha_i', \alpha_k}$
is $-1-h(\square_i)$; the weight of $f_{\alpha_k', \alpha_k}$ is
$1$; if $(k+1) \le i \le s$, then the weight of $f_{\alpha_i',
\alpha_k}$ is $1+h(\square_i)$. Note that the weight of
$\phi(f_{\alpha_i, \beta_i})$ is the same as the weight of
$f_{\alpha_i, \beta_i}$. Hence if $0 \le i \le (k-2)$, then the
weight of $\phi(f_{\alpha_i, \beta_i})$ is $(p_i-1)-h(\square_i) =
-h(\square_i')$; if $k \le i \le s-1$, then the weight of
$\phi(f_{\alpha_{i+1}, \beta_{i+1}})$ is
$-(q_{i+1}-1)+h(\square_{i+1}) = h(\square_{i+1}')$. By
(\ref{case1a_ker.1}) and (\ref{hil_ts.2}),
\begin{eqnarray}    \label{case1a_ker.2}
    e_{\mathbb T}
&=&(-1)^{n+1} \cdot h(\la)^2
   \cdot \prod_{0 \le i \le m \atop i \ne k-1, k}
   \frac{1+h(\square_i)}{h(\square_i')} \cdot
   [1+h(\square_{k-1})] \cdot t^{2(n+1)}   \nonumber \\
&=&(-1)^{n+1} \cdot h(\la)^2
   \cdot \prod_{0 \le i \le m \atop i \ne k}
   \frac{1+h(\square_i)}{h(\square_i')} \cdot t^{2(n+1)}
\end{eqnarray}
where we have used the observation that $h(\square_{k-1}') = 1$.

\medskip\noindent
{\bf Case 1b:} $p_k=1$ but $q_k \ne 1$. Then $s=m$, $\alpha_k'
\not\in A$, and $\alpha_i' = \alpha_i$ if $0 \le i \le s$ and $i \ne
k$. For $0 \le i \le (k-1)$, let $\beta_i \in B$ be the $p_i$-th
cell directly above $\alpha_k$. Then, $\beta_i \in P_{\alpha_i}$.
For $(k+1) \le i \le (s-1)$, let $\beta_{i+1} \in B$ be the
$q_{i+1}$-th cell directly to the left of $\alpha_k$. Then,
$\beta_{i+1} \in Q_{\alpha_{i+1}}$. As in Case 1a, there is an exact
sequence
\begin{eqnarray}    \label{case1b_ker.1}
0 \to \frac{\ker(\phi-\psi)}{\bigoplus_{i=0}^s
  \C f_{\alpha_i', \alpha_k}} \to \Hom \big (I_{\xi_{\la}},
  R/I_{\xi_{\la}} \big ) \oplus {\rm im}(\psi)
  \qquad\qquad   \nonumber   \\
\to {\rm im}(\psi) \bigoplus \left ( \bigoplus_{i = 0}^{k-1}
  \C \phi(f_{\alpha_i, \beta_i}) \right )
  \bigoplus \left ( \bigoplus_{i = k+1}^{s-1}
  \C \phi(f_{\alpha_{i+1}, \beta_{i+1}}) \right ) \to 0.
\end{eqnarray}

If $0 \le i \le (k-1)$, then the weight of $f_{\alpha_i', \alpha_k}$
is $-1-h(\square_i)$; the weight of $f_{\alpha_k', \alpha_k}$ is
$-1$; if $(k+1) \le i \le s$, then the weight of $f_{\alpha_i',
\alpha_k}$ is $1+h(\square_i)$. If $0 \le i \le (k-1)$, then the
weight of $\phi(f_{\alpha_i, \beta_i})$ is $-h(\square_i')$; if $k+1
\le i \le s-1$, then the weight of $\phi(f_{\alpha_{i+1},
\beta_{i+1}})$ is $h(\square_{i+1}')$. Combining with
$h(\square_{k+1}') = 1$, we obtain
\begin{eqnarray}    \label{case1b_ker.2}
    e_{\mathbb T}
&=&(-1)^{n+1} \cdot h(\la)^2
   \cdot \prod_{0 \le i \le m \atop i \ne k, k+1}
   \frac{1+h(\square_i)}{h(\square_i')} \cdot
   [1+h(\square_{k+1})] \cdot t^{2(n+1)}  \nonumber \\
&=&(-1)^{n+1} \cdot h(\la)^2
   \cdot \prod_{0 \le i \le m \atop i \ne k}
   \frac{1+h(\square_i)}{h(\square_i')} \cdot t^{2(n+1)}.
\end{eqnarray}

\medskip\noindent
{\bf Case 2:} $k = 0$ and $p_0 > 1$, or $k=m$ and $q_m > 1$, or $0 <
k < m$ and $p_k, q_k > 1$. Then $s=m+1$, $\alpha_k', \alpha_{k+1}'
\not\in A$, $\alpha_i' = \alpha_i$ if $0 \le i \le (k-1)$, and
$\alpha_i' = \alpha_{i-1}$ if $(k+2) \le i \le s$. For $0 \le i \le
(k-1)$, we have $p_i' = p_i$ and let $\beta_i \in B$ be the $p_i$-th
cell directly above $\alpha_k$. Then, $\beta_i \in P_{\alpha_i}$.
For $(k+1) \le i \le (s-1)$, we have $q_{i+1}' = q_{i}$ and
$\alpha_{i+1}' = \alpha_{i}$. Let $\beta_{i} \in B$ be the
$q_{i}$-th cell directly to the left of $\alpha_k$. Then, $\beta_{i}
\in Q_{\alpha_{i}}$. There is an exact sequence
\begin{eqnarray}    \label{case2_ker.1}
0 \to \frac{\ker(\phi-\psi)}{\bigoplus_{i=0}^s
  \C f_{\alpha_i', \alpha_k}} \to \Hom \big (I_{\xi_{\la}},
  R/I_{\xi_{\la}} \big ) \oplus {\rm im}(\psi)
  \qquad\qquad   \nonumber   \\
\to {\rm im}(\psi) \bigoplus \left ( \bigoplus_{i = 0}^{k-1}
  \C \phi(f_{\alpha_i, \beta_i}) \right )
  \bigoplus \left ( \bigoplus_{i = k+1}^{s-1}
  \C \phi(f_{\alpha_{i}, \beta_{i}}) \right ) \to 0.
\end{eqnarray}

If $0 \le i \le (k-1)$, then the weight of $f_{\alpha_i', \alpha_k}$
is $-1-h(\square_i)$; the weight of $f_{\alpha_k', \alpha_k}$ is
$-1$; the weight of $f_{\alpha_{k+1}', \alpha_k}$ is $1$; if $(k+2)
\le i \le s$, then the weight of $f_{\alpha_i', \alpha_k}$ is
$1+h(\square_{i-1})$. If $0 \le i \le (k-1)$, then the weight of
$\phi(f_{\alpha_i, \beta_i})$ is $-h(\square_i')$; if $k+1 \le i \le
s-1$, then the weight of $\phi(f_{\alpha_i, \beta_i})$ is
$h(\square_i')$. Hence (\ref{lemma_eT.1}) holds.

\medskip\noindent
{\bf Case 3:} $p_k= q_k = 1$. Then, $s = m-1$, $\alpha_i' =
\alpha_i$ if $0 \le i \le (k-1)$, and $\alpha_i' = \alpha_{i+1}$ if
$k \le i \le s$. For $0 \le i \le (k-2)$, we have $p_i' = p_i$ and
let $\beta_i \in B$ be the $p_i$-th cell directly above $\alpha_k$.
Then, $\beta_i \in P_{\alpha_i}$. For $k \le i \le (s-1)$, we have
$q_{i+1}' = q_{i+2}$ and $\alpha_{i+1}' = \alpha_{i+2}$. Let
$\beta_{i+2} \in B$ be the $q_{i+2}$-th cell directly to the left of
$\alpha_k$. Then, $\beta_{i+2} \in Q_{\alpha_{i+2}}$. There is an
exact sequence
\begin{eqnarray}    \label{case3_ker}
0 \to \frac{\ker(\phi-\psi)}{\bigoplus_{i=0}^s
  \C f_{\alpha_i', \alpha_k}} \to \Hom \big (I_{\xi_{\la}},
  R/I_{\xi_{\la}} \big ) \oplus {\rm im}(\psi)
  \qquad\qquad   \nonumber   \\
\to {\rm im}(\psi) \bigoplus \left ( \bigoplus_{i = 0}^{k-2}
  \C \phi(f_{\alpha_i, \beta_i}) \right )
  \bigoplus \left ( \bigoplus_{i = k}^{s-1}
  \C \phi(f_{\alpha_{i+2}, \beta_{i+2}}) \right ) \to 0.
\end{eqnarray}

If $0 \le i \le (k-1)$, then the weight of $f_{\alpha_i', \alpha_k}$
is $-1-h(\square_i)$; if $k \le i \le s$, then the weight of
$f_{\alpha_i', \alpha_k}$ is $1+h(\square_{i+1})$. If $0 \le i \le
(k-2)$, then the weight of $\phi(f_{\alpha_i, \beta_i})$ is
$-h(\square_i')$; if $k \le i \le s-1$, then the weight of
$\phi(f_{\alpha_{i+2}, \beta_{i+2}})$ is $h(\square_{i+2}')$. Hence
\begin{eqnarray*}
(-1)^{n+1} \cdot h(\la)^2 \cdot \prod_{0 \le i \le m \atop i \ne
k-1, k, k+1} \frac{1+h(\square_i)}{h(\square_i')} \cdot
[1+h(\square_{k-1})] \cdot [1+h(\square_{k+1})] \cdot t^{2(n+1)}
\end{eqnarray*}
noting that $h(\square_{k-1}') = h(\square_{k+1}') = 1$. Therefore,
(\ref{lemma_eT.1}) holds.



\end{document}